# Colimits of Enriched Categories

Andy Sukowski-Bang

Abstract. We introduce categories of strings to give an explicit hom-object formula for colimits of enriched categories. As corollaries, we recover the necklace construction by Dugger and Spivak and describe colimits in Cat and of Lawvere metric spaces.

## Contents

## 1. Introduction

Let a monoidal category $\mathcal{V}$ be cocomplete and let the tensor product $\otimes$ be cocontinuous in each argument. This paper explicitly constructs the colimit of a diagram $F\colon \mathcal{I} \to \mathrm{Cat}_{\mathcal{V}}$ in the category of $\mathcal{V}$-enriched categories in Theorem 4.14. [Wol74] proves the existence of $\mathrm{colim}(F)$, but gives no explicit formula for hom-objects $(\mathrm{colim}\, F)(a, b)$. In the special case of a certain diagram $F\colon (\Delta \downarrow X) \to \mathrm{Cat}_{\mathrm{sSet}}$ of simplicially enriched categories, Dugger and Spivak provide an explicit construction via necklaces in [DS11]. Section 6 recovers this special case in Corollary 6.18 [DS11, Corollary 4.4].

> In [Lur08, p. 29] each $\mathfrak{C}X$ is defined as a certain colimit in the category $\mathrm{Cat}_{\mathrm{sSet}}$, but colimits in $\mathrm{Cat}_{\mathrm{sSet}}$ are notoriously difficult to understand.
>
> — D. Dugger and D. Spivak [DS11, §1.1]

Before deriving the necklace construction, we warm up on metric spaces and ordinary categories in Section 5. To build intuition for colimits of diagrams $\mathcal{I} \to \mathrm{Cat}_{\mathcal{V}}$, we now look at a coequalizer and a pushout in ordinary Cat.

Example 1.1. Consider the category $\mathcal{I}$ consisting of two morphisms $f, g\colon i \to j$ and the following diagram $F\colon \mathcal{I} \to \mathrm{Cat}$ where $y = (Ff)(x)$ and $z = (Fg)(x)$, and the indicated morphisms generate the categories. Then the coequalizer $\mathrm{colim}(F)$ with colimiting cocone $\iota\colon F \Rightarrow \Delta_{\mathcal{I}}(\mathrm{colim}\, F)$ has a single object $a = \iota_i(x) = \iota_j(y) = \iota_j(z)$ but infinitely many endomorphisms of $a$, namely the compositions $\iota_j(h) \circ \cdots \circ \iota_j(h)$.

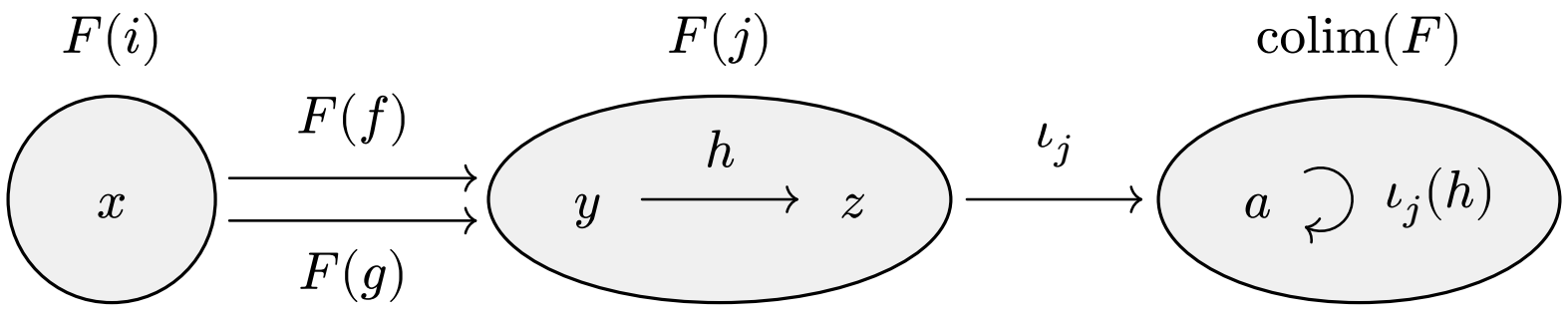

Example 1.2. Let $\mathcal{I}$ be the span $k \overset{g}{\leftarrow} i \overset{f}{\rightarrow} j$ and $F\colon \mathcal{I} \to \mathrm{Cat}$ the diagram below, where dashed arrows indicate compositions. The pushout $\mathrm{colim}(F)$ with colimiting cocone $\iota\colon F \Rightarrow \Delta_{\mathcal{I}}(\mathrm{colim}\, F)$ has objects $a = \iota_i(x)$, $b = \iota_i(y)$ and $c = \iota_i(z)$. The hom-sets $(Fi)(x, y)$ and $(Fj)((Ff)(x), (Ff)(y))$ are "glued together". Note that we obtain infinitely many compositions, freely generated by the indicated morphisms in $\mathrm{colim}(F)$.

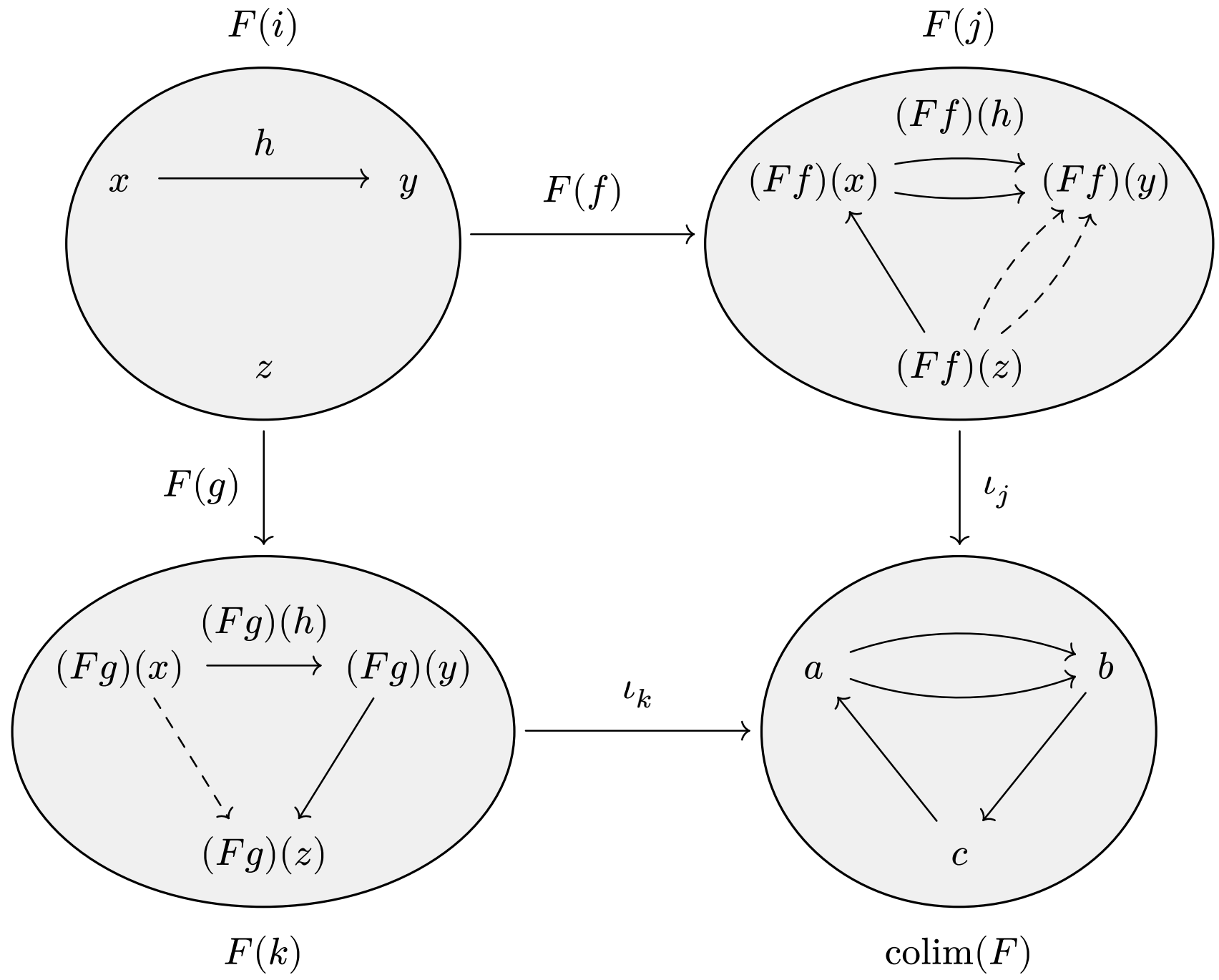

Observation 1.3. Hom-objects become composable if endpoints coincide in $\mathrm{colim}(\mathrm{Ob} \circ F)$. In Example 1.2, the hom-objects of $F(i)$ and $F(j)$ are "glued together" by a transport along the underlying morphism $f$ of $\mathcal{I}$. We cannot just take the colimit of graphs followed by free compositions, since old compositions in the individual $F(i)$ need to be preserved.

For a general diagram $F\colon \mathcal{I} \to \mathrm{Cat}_{\mathcal{V}}$, we obtain objects $\mathrm{Ob}(\mathrm{colim}\, F) = \mathrm{colim}(\mathrm{Ob} \circ F)$. In Section 4, we construct the hom-object $(\mathrm{colim}\, F)(a, b)$ as a colimit in $\mathcal{V}$ of tensored hom-objects

$$F(i_n)(x_n, y_n) \otimes \cdots \otimes F(i_1)(x_1, y_1), \tag{1}$$

specifically, the colimit of the **string evaluation** diagram $\mathrm{ev}^F_{a,b}\colon \mathcal{S}_F(a, b) \to \mathcal{V}$ (see Definition 4.8) over the **category of strings** $\mathcal{S}_F(a, b)$ (see Definition 4.5). A non-empty $F$-string $s \in \mathcal{S}_F(a, b)$ is of the form

$$s = [(i_n, x_n, y_n), ..., (i_1, x_1, y_1)]$$

between $a = [x_1] \in \mathrm{colim}(\mathrm{Ob} \circ F)$ and $b = [y_n] \in \mathrm{colim}(\mathrm{Ob} \circ F)$, with $x_k, y_k \in F(i_k)$ and compatible endpoints $[y_{k-1}] = [x_k] \in \mathrm{colim}(\mathrm{Ob} \circ F)$. Its string evaluation $\mathrm{ev}^F_{a,b}(s)$ is exactly Equation 1. The question remains: how to define composition $\circ_{a,b,c}$ in $\mathrm{colim}(F)$? Abbreviating $\mathcal{S} := \mathcal{S}_F$ and suppressing endpoints in $\mathrm{ev} := \mathrm{ev}^F$, we get the first isomorphism from cocontinuity of $\otimes$ in each argument. But how do we obtain $\varphi$?

$$\operatorname*{colim}_{\mathcal{S}(b,c)} \mathrm{ev} \otimes \operatorname*{colim}_{\mathcal{S}(a,b)} \mathrm{ev} \cong \operatorname*{colim}_{\mathcal{S}(b,c)\times\mathcal{S}(a,b)} (\mathrm{ev} \boxtimes \mathrm{ev}) \overset{\varphi}{\longrightarrow} \operatorname*{colim}_{\mathcal{S}(a,c)} \mathrm{ev}.$$

The domain of $\varphi$ is a colimit of a diagram of shape $\mathcal{S}(b, c) \times \mathcal{S}(a, b)$, while the codomain's diagram has shape $\mathcal{S}(a, c)$. In Section 3 we will show that taking colimits is a strong monoidal functor on the Grothendieck construction of $[-, \mathcal{V}]$, keeping track of shape changes $\mathcal{J} \to \mathcal{J}'$.

$$\mathrm{colim}\colon \int [-, \mathcal{V}] \to \mathcal{V}, \quad (\mathcal{J}, G) \mapsto \operatorname*{colim}_{\mathcal{J}} G$$

To prove the universal property of our explicit colimit construction, we will encode the unit maps $\mathbb{1} \to (\mathrm{colim}\, F)(a, a)$ directly in $\mathcal{S}_F$. For applications, we introduce the smaller **category of non-empty strings** $\overline{\mathcal{S}}_F$ and define unit maps externally in Section 4.3.

## 2. Enriched categories and change of base

A category enriched in $\mathcal{V}$ has hom-objects $\mathcal{C}(a,b) \in \mathcal{V}$ instead of ordinary hom-sets. That leaves the question of how hom-objects are composed; for that, we need a tensor product $\otimes$. The tensor unit $\mathbb{1}$ is used to define identities.

Definition 2.1 (monoidal category). A category $\mathcal{V}$ is **monoidal** when equipped with

- a tensor product functor $\otimes : \mathcal{V} \times \mathcal{V} \to \mathcal{V}$,
- a tensor unit $\mathbb{1} \in \mathcal{V}$,
- an associator, i.e. a natural isomorphism $\alpha_{a,b,c} \colon (a \otimes b) \otimes c \to a \otimes (b \otimes c)$,
- a left unitor, i.e. a natural isomorphism $l_a \colon \mathbb{1} \otimes a \to a$,
- a right unitor, i.e. a natural isomorphism $r_a \colon a \otimes \mathbb{1} \to a$,

satisfying the triangle and pentagon identities.

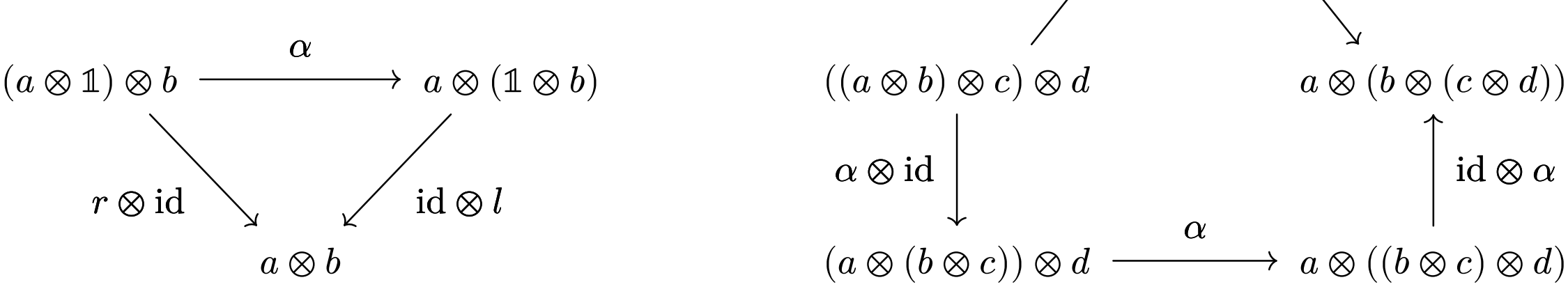

To obtain a $\mathcal{V}$-enriched category, we need to define composition morphisms, as well as identities. Since $\mathcal{V}$ may not necessarily have an underlying set structure, we cannot simply choose an identity morphism $\mathrm{id}_a$ from a hom-object $\mathcal{C}(a,a) \in \mathcal{V}$. Instead, we pick an identity via a morphism $j_a \colon \mathbb{1} \to \mathcal{C}(a,a)$ in $\mathcal{V}$.

Definition 2.2 (enriched category). A small category $\mathcal{C}$ **enriched** in a monoidal category $\mathcal{V}$ is

- a set $\mathrm{Ob}(\mathcal{C})$ often denoted by $\mathcal{C}$,
- for each tuple $(a,b) \in \mathcal{C} \times \mathcal{C}$ a hom-object $\mathcal{C}(a,b) \in \mathcal{V}$,
- for each triple $(a,b,c)$ a composition morphism $\circ_{a,b,c} : \mathcal{C}(b,c) \otimes \mathcal{C}(a,b) \to \mathcal{C}(a,c)$ in $\mathcal{V}$,
- for each $a \in \mathcal{C}$ a morphism $j_a \colon \mathbb{1} \to \mathcal{C}(a,a)$ called the unit map,

with commuting associativity and unital diagrams.

$$\begin{array}{ccc} (\mathcal{C}(c,d) \otimes \mathcal{C}(b,c)) \otimes \mathcal{C}(a,b) & \xrightarrow{\alpha} & \mathcal{C}(c,d) \otimes (\mathcal{C}(b,c) \otimes \mathcal{C}(a,b)) \\ \downarrow {\scriptstyle \circ \otimes \mathrm{id}} & & \downarrow {\scriptstyle \mathrm{id} \otimes \circ} \\ \mathcal{C}(b,d) \otimes \mathcal{C}(a,b) & & \mathcal{C}(c,d) \otimes \mathcal{C}(a,c) \\ {\scriptstyle \circ} \searrow & & \swarrow {\scriptstyle \circ} \\ & \mathcal{C}(a,d) & \end{array}$$

$$\begin{array}{ccc} \mathbb{1} \otimes \mathcal{C}(a,b) & \xrightarrow{j_b \otimes \mathrm{id}} & \mathcal{C}(b,b) \otimes \mathcal{C}(a,b) \\ {\scriptstyle l} \searrow & & \swarrow {\scriptstyle \circ} \\ & \mathcal{C}(a,b) & \end{array}$$

$$\begin{array}{ccc} \mathcal{C}(a,b) \otimes \mathbb{1} & \xrightarrow{\mathrm{id} \otimes j_a} & \mathcal{C}(a,b) \otimes \mathcal{C}(a,a) \\ {\scriptstyle r} \searrow & & \swarrow {\scriptstyle \circ} \\ & \mathcal{C}(a,b) & \end{array}$$

To define the category of small $\mathcal{V}$-enriched categories, we introduce the notion of an enriched functor.

Definition 2.3 (enriched functor). Let $\mathcal{C}$ and $\mathcal{D}$ be small $\mathcal{V}$-enriched categories. An **enriched functor** $F \colon \mathcal{C} \to \mathcal{D}$ consists of a function $F \colon \mathrm{Ob}(\mathcal{C}) \to \mathrm{Ob}(\mathcal{D})$, and a collection of morphisms $F_{a,b} \colon \mathcal{C}(a,b) \to \mathcal{D}(Fa, Fb)$, respecting composition and units.

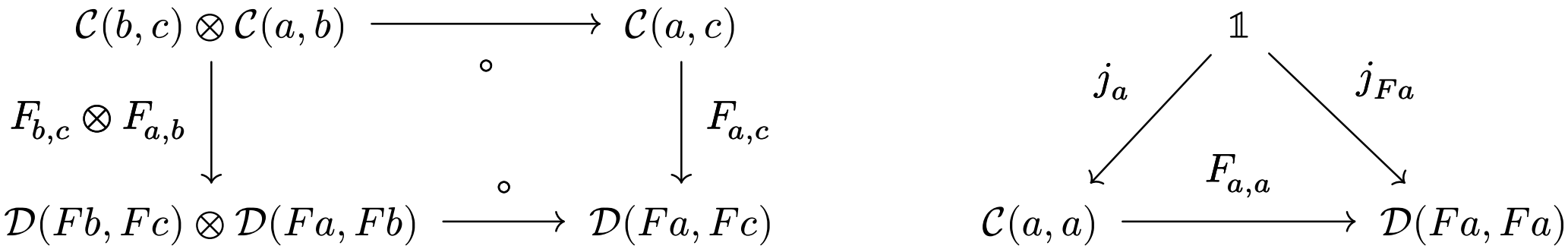

One can change the enriching base via a lax monoidal functor $F \colon \mathcal{U} \to \mathcal{V}$, which induces a functor $F_* \colon \mathrm{Cat}_{\mathcal{U}} \to \mathrm{Cat}_{\mathcal{V}}$.

Definition 2.4 (lax monoidal functor). Let $(\mathcal{U}, \times, *)$ and $(\mathcal{V}, \otimes, \mathbb{1})$ be monoidal categories. A functor $F\colon \mathcal{U} \to \mathcal{V}$ is **lax monoidal** if equipped with coherence maps:

- a morphism $\eta\colon \mathbb{1} \to F(*)$,
- a natural transformation $\mu_{a,b}\colon F(a) \otimes F(b) \to F(a \times b)$,

satisfying unitality and associativity.

$$\begin{CD} \mathbb{1}\otimes Fa @>{\eta\otimes \mathrm{id}}>> F(*)\otimes Fa \\ @V{l}VV @VV{\mu}V \\ Fa @<{Fl}<< F(*\times a)\end{CD} \qquad \begin{CD} Fa\otimes \mathbb{1} @>{\mathrm{id}\otimes\eta}>> Fa\otimes F(*) \\ @V{r}VV @VV{\mu}V \\ Fa @<{Fr}<< F(a\times *)\end{CD}$$

$$\begin{CD} (Fa\otimes Fb)\otimes Fc @>{\alpha}>> Fa\otimes(Fb\otimes Fc) \\ @V{\mu\otimes\mathrm{id}}VV @VV{\mathrm{id}\otimes\mu}V \\ F(a\times b)\otimes Fc @. Fa\otimes F(b\times c) \\ @V{\mu}VV @VV{\mu}V \\ F((a\times b)\times c) @>{F\alpha}>> F(a\times(b\times c))\end{CD}$$

If $\eta$ and $\mu$ are isomorphisms, $F$ is called **strong monoidal**.

In Theorem 4.14, our construction of $\operatorname{colim}(F)$ for a diagram $F\colon \mathcal{I} \to \mathrm{Cat}_{\mathcal{V}}$ emerges by change of enriching base along a strong monoidal functor $\int[-, \mathcal{V}] \to \mathcal{V}$ defined in Proposition 3.8. In Definition 6.7, the rigidification of a standard simplex is defined via a change of base $\mathrm{Cat} \to \mathrm{sSet}$.

Lemma 2.5 (change of enriching base). A lax monoidal functor $F\colon \mathcal{U} \to \mathcal{V}$ induces a functor $F_*\colon \mathrm{Cat}_{\mathcal{U}} \to \mathrm{Cat}_{\mathcal{V}}$, sending a $\mathcal{U}$-enriched category $\mathcal{C}$ to $F_*(\mathcal{C})$ with $\mathrm{Ob}(F_*\mathcal{C}) \coloneqq \mathrm{Ob}(\mathcal{C})$ and $F_*\mathcal{C}(a,b) \coloneqq F(\mathcal{C}(a,b))$, as well as composition and unit maps

$$F(\mathcal{C}(b,c)) \otimes F(\mathcal{C}(a,b)) \xrightarrow{\mu} F(\mathcal{C}(b,c) \times \mathcal{C}(a,b)) \xrightarrow{F(\circ)} F(\mathcal{C}(a,c)), \qquad \mathbb{1} \xrightarrow{\eta} F(*) \xrightarrow{F(j_a)} F(\mathcal{C}(a,a)).$$

Let $G\colon \mathcal{C} \to \mathcal{D}$ be a $\mathcal{U}$-enriched functor. Then $(F_*G)_{a,b} \coloneqq F(G_{a,b})$ defines a $\mathcal{V}$-enriched functor.

Remark 2.6. See [Cru08, §4.2] for 2-functoriality of $F_*$.

Proof of Lemma 2.5. We will only prove that $F_*(\mathcal{C})$ satisfies associativity and unitality. By naturality, $F(\circ \times \mathrm{id}) \circ \mu = \mu \circ (F(\circ) \otimes F(\mathrm{id})) = \mu \circ (F(\circ) \otimes \mathrm{id})$.

$$\begin{CD} (F\mathcal{C}(c,d)\otimes F\mathcal{C}(b,c))\otimes F\mathcal{C}(a,b) @>>{\alpha}> F\mathcal{C}(c,d)\otimes(F\mathcal{C}(b,c)\otimes F\mathcal{C}(a,b)) \\ @V{\mu\otimes\mathrm{id}}VV @VV{\mathrm{id}\otimes\mu}V \\ F(\mathcal{C}(c,d)\times\mathcal{C}(b,c))\otimes F\mathcal{C}(a,b) @. F\mathcal{C}(c,d)\otimes F(\mathcal{C}(b,c)\times\mathcal{C}(a,b)) \\ @V{\mu}VV @VV{\mu}V \\ F((\mathcal{C}(c,d)\times\mathcal{C}(b,c))\times\mathcal{C}(a,b)) @>{F\alpha}>> F(\mathcal{C}(c,d)\times(\mathcal{C}(b,c)\times\mathcal{C}(a,b))) \\ @V{F(\circ\times\mathrm{id})}VV @VV{F(\mathrm{id}\times\circ)}V \\ F(\mathcal{C}(b,d)\times\mathcal{C}(a,b)) @. F(\mathcal{C}(c,d)\times\mathcal{C}(a,c)) \end{CD}$$

$$F(\mathcal{C}(b,d)\times\mathcal{C}(a,b)) \xrightarrow{F(\circ)} F\mathcal{C}(a,d) \xleftarrow{F(\circ)} F(\mathcal{C}(c,d)\times\mathcal{C}(a,c))$$

By naturality, $F(\mathrm{id} \times j_a) \circ \mu = \mu \circ (F(\mathrm{id}) \otimes F(j_a)) = \mu \circ (\mathrm{id} \otimes F(j_a))$.

$$\begin{CD} F\mathcal{C}(a,b)\otimes\mathbb{1} @>>{\mathrm{id}\otimes\eta}> F\mathcal{C}(a,b)\otimes F(*) \\ @V{r}VV @VV{\mu}V \\ F\mathcal{C}(a,b) @<{Fr}<< F(\mathcal{C}(a,b)\times *) \end{CD}$$

$$F\mathcal{C}(a,b) \xleftarrow{F(\circ)} F(\mathcal{C}(a,b)\times\mathcal{C}(a,a)) \xleftarrow{F(\mathrm{id}\times j_a)} F(\mathcal{C}(a,b)\times *)$$

□

## 3. Monoidal colimit functoriality

Recall from the introduction that we want to define composition in Theorem 4.14 via cocontinuity of $\otimes$ in each argument. Our problem is the indicated morphism $\varphi$.

$$\operatorname*{colim}_{\mathcal{S}(b,c)} \mathrm{ev} \otimes \operatorname*{colim}_{\mathcal{S}(a,b)} \mathrm{ev} \cong \operatorname*{colim}_{\mathcal{S}(b,c)\times\mathcal{S}(a,b)} (\mathrm{ev} \boxtimes \mathrm{ev}) \xrightarrow{\varphi} \operatorname*{colim}_{\mathcal{S}(a,c)} \mathrm{ev}$$

The domain of $\varphi$ is a colimit of a diagram of shape $\mathcal{S}(b,c) \times \mathcal{S}(a,b)$, while the codomain's diagram has shape $\mathcal{S}(a,c)$. We will see that taking colimits is functorial on the Grothendieck construction of $[-,\mathcal{V}]$, keeping track of shape changes $\mathcal{J} \to \mathcal{J}'$. In fact, we will even show strong monoidal functoriality in Proposition 3.8, enabling change of enriching base $\mathrm{Cat}_{\int[-,\mathcal{V}]} \to \mathrm{Cat}_{\mathcal{V}}$ (see Lemma 2.5).

Definition 3.1 (Grothendieck construction) [GR71, §VI.8]. Let $F\colon \mathcal{C}^{\mathsf{op}} \to \mathrm{CAT}$ be a functor. The (contravariant) **Grothendieck construction** $\int F$ is the category whose
- objects are pairs $(c, e)$ of $c \in \mathcal{C}$ and $e \in F(c)$,
- morphisms $(c, e) \to (c', e')$ are given by pairs $(f\colon c \to c', e \to (Ff)(e'))$,
- compositions are $(g, \psi) \circ (f, \varphi) := (gf, (Ff)(\psi) \circ \varphi)$,
- identities are $(\mathrm{id}_c, \mathrm{id}_e)$.

Remark 3.2. We allow $F$ to have values in CAT, the category of large categories, since we are interested in the Grothendieck construction of $[-,\mathcal{V}]\colon \mathrm{Cat}^{\mathsf{op}} \to \mathrm{CAT}$. Even if $\mathcal{I}$ is small, the functor category $[\mathcal{I},\mathcal{V}]$ may be large.

Proposition 3.3 (functoriality of colimit). Let $\mathcal{V}$ be cocomplete. A choice of colimit per diagram with small shape category defines a functor

$$\mathrm{colim}\colon \int [-,\mathcal{V}] \to \mathcal{V}, \quad (\mathcal{I}, F) \mapsto \mathrm{colim}\, F$$

sending $H\colon \mathcal{I} \to \mathcal{J}$ and $\eta\colon F \Rightarrow G \circ H$ to a unique map $\mathrm{colim}(H,\eta)\colon \mathrm{colim}\, F \to \mathrm{colim}\, G$.

Remark 3.4. This generalizes functoriality of colim: $[\mathcal{I},\mathcal{V}] \to \mathcal{V}$ for any small category $\mathcal{I}$ [Rie16, Proposition 3.3.1] by precomposing $[\mathcal{I},\mathcal{V}] \hookrightarrow \int[-,\mathcal{V}]$, sending $F \mapsto (\mathcal{I}, F)$ and $\eta \mapsto (\mathrm{id}, \eta)$.

Proof of Proposition 3.3. $(\mathcal{I}, F) \to (\mathcal{J}, G)$ consists of a functor $H\colon \mathcal{I} \to \mathcal{J}$ and a natural transformation $\eta\colon F \Rightarrow [H,\mathcal{V}](G) = G \circ H$. Whiskering the colimiting cocone $\iota\colon G \Rightarrow \Delta_{\mathcal{J}} \operatorname{colim}(G)$ with $H$, and precomposing with $\eta$, we obtain the cocone

$$F \overset{\eta}{\Longrightarrow} G \circ H \overset{\iota H}{\Longrightarrow} (\Delta_{\mathcal{J}} \operatorname{colim} G) \circ H = \Delta_{\mathcal{I}} \operatorname{colim} G.$$

The universal property of $\operatorname{colim}(F)$ implies a unique cocone morphism $\operatorname{colim}(H,\eta)\colon \operatorname{colim} F \to \operatorname{colim} G$. Consider the composition

$$(\mathcal{I}, F) \xrightarrow{(H,\eta)} (\mathcal{I}', F') \xrightarrow{(H',\eta')} (\mathcal{I}'', F'')$$

and let $\iota$, $\iota'$ and $\iota''$ denote the respective colimits' cocones. By definition, $\Delta \operatorname{colim}(H,\eta) \circ \iota = \iota' H \circ \eta$ and $\Delta \operatorname{colim}(H',\eta') \circ \iota' = \iota'' H' \circ \eta'$. Applying colim to the composition in $\int[-,\mathcal{V}]$, we have:

$$\begin{aligned} \Delta \operatorname{colim}(H'H, [H,\mathcal{V}](\eta') \circ \eta) \circ \iota &= \Delta \operatorname{colim}(H'H, \eta' H \circ \eta) \circ \iota = \iota'' H' H \circ \eta' H \circ \eta \\ = (\iota'' H' \circ \eta') H \circ \eta &= \Delta \operatorname{colim}(H',\eta') \circ \iota' H \circ \eta = \Delta \operatorname{colim}(H',\eta') \circ \Delta \operatorname{colim}(H,\eta) \circ \iota \end{aligned}$$

That is commutativity of the following diagram in $[\mathcal{I},\mathcal{V}]$.

$$\begin{array}{ccccc} F & \overset{\eta}{\Longrightarrow} & F'H & \overset{\eta' H}{\Longrightarrow} & F''H'H \\ \Downarrow \iota & & \Downarrow \iota' H & & \Downarrow \iota'' H' H \\ \Delta \operatorname{colim}_{\mathcal{I}} F & \overset{\Delta \operatorname{colim}(H,\eta)}{\Longrightarrow} & \Delta \operatorname{colim}_{\mathcal{I}'} F' & \overset{\Delta \operatorname{colim}(H',\eta')}{\Longrightarrow} & \Delta \operatorname{colim}_{\mathcal{I}''} F'' \end{array}$$

□

Lemma 3.5. A monoidal category $(\mathcal{V}, \otimes, \mathbb{1})$ canonically endows $\int[-, \mathcal{V}]$ with a monoidal structure. Let $F \boxtimes G := \otimes \circ (F \times G)$. Overloading notation, define the tensor product $(\mathcal{I}, F) \boxtimes (\mathcal{J}, G) := (\mathcal{I} \times \mathcal{J}, F \boxtimes G)$ with the tensor unit $(*, * \mapsto \mathbb{1})$ on the terminal category.

Proof. Use the cartesian monoidal structure of $(\mathrm{Cat}, \times, *)$. The associator in $\int[-, \mathcal{V}]$ has components $\left(\alpha^{\mathrm{Cat}}_{\mathcal{I},\mathcal{J},\mathcal{K}}, \alpha^{\mathcal{V}}_{F,G,H}\right)$, the right unitor has components $\left(r^{\mathrm{Cat}}_{\mathcal{I}}, r^{\mathcal{V}}_{F}\right)$ and the left unitor has components $\left(l^{\mathrm{Cat}}_{\mathcal{I}}, l^{\mathcal{V}}_{F}\right)$. Note that we are whiskering $\alpha^{\mathcal{V}}_{F,G,H} = \alpha^{\mathcal{V}}((F \times G) \times H)\colon (F \boxtimes G) \boxtimes H \Rightarrow (F \boxtimes (G \boxtimes H)) \circ \alpha^{\mathrm{Cat}}_{\mathcal{I},\mathcal{J},\mathcal{K}}$. We verify the triangle and pentagon identities in Lemma A.1 of Appendix A ◻

Definition 3.6 (cocontinuity). A functor $G\colon \mathcal{C} \to \mathcal{D}$ is **cocontinuous** if for any small diagram $F\colon \mathcal{I} \to \mathcal{C}$ and its colimiting cocone $\kappa\colon F \Rightarrow \Delta\,\mathrm{colim}(F)$, the induced morphism of the cocone

$$G \circ F \overset{G\kappa}{\Longrightarrow} G \circ \Delta\,\mathrm{colim}\,F = \Delta G(\mathrm{colim}\,F)$$

is an isomorphism $\mathrm{colim}(G \circ F) \xrightarrow{\cong} G(\mathrm{colim}\,F)$.

Example 3.7. Let $G\colon \mathcal{C} \times \mathcal{C}' \to \mathcal{D}$ be cocontinuous in each argument. For two small diagrams $F\colon \mathcal{I} \to \mathcal{C}$ and $F'\colon \mathcal{I}' \to \mathcal{C}'$ with colimiting cocones $\kappa$ and $\kappa'$, the cocone

$$G \circ (F \times F') \overset{G(\kappa\times\kappa')}{\Longrightarrow} G \circ (\Delta_{\mathcal{I}}\,\mathrm{colim}\,F \times \Delta_{\mathcal{I}'}\,\mathrm{colim}\,F') = \Delta_{\mathcal{I}\times\mathcal{I}'} G(\mathrm{colim}\,F, \mathrm{colim}\,F')$$

induces a morphism $\mathrm{colim}(G \circ (F \times F')) \to G(\mathrm{colim}\,F, \mathrm{colim}\,F')$. This is an isomorphism, since we can apply Fubini for colimits [Mac98, §IX.2] to the bifunctor $G(F, F') := G \circ (F \times F')$, followed by cocontinuity in $\mathcal{C}'$, then $\mathcal{C}$.

$$\underset{\mathcal{I}\times\mathcal{I}'}{\mathrm{colim}}\, G(F, F') \cong \underset{\mathcal{I}}{\mathrm{colim}}\,\underset{\mathcal{I}'}{\mathrm{colim}}\, G(F, F') \cong \underset{\mathcal{I}}{\mathrm{colim}}\, G(F, \mathrm{colim}\,F') \cong G(\mathrm{colim}\,F, \mathrm{colim}\,F')$$

Proposition 3.8. Let $\mathcal{V}$ be cocomplete and $\otimes$ be cocontinuous in each argument. The colimit functor $\mathrm{colim}\colon \int[-, \mathcal{V}] \to \mathcal{V}$ from Proposition 3.3 extends to a strong monoidal functor with a natural isomorphism $\mu_{F,G}\colon \mathrm{colim}\,F \otimes \mathrm{colim}\,G \cong \mathrm{colim}(F \boxtimes G)$ and an isomorphism $\eta\colon \mathbb{1} \cong \mathrm{colim}(* \mapsto \mathbb{1})$.

Proof. We prove unital and associative coherence in Lemma A.2 of Appendix A. Take $\eta$ to be the only component in the cocone of $\mathrm{colim}(* \mapsto \mathbb{1})$. By cocontinuity of $\otimes$, the component $\mu^{-1}_{F,G}$ induced by the cocone

$$F \boxtimes F' \overset{\kappa\boxtimes\kappa'}{\Longrightarrow} (\Delta_{\mathcal{I}}\,\mathrm{colim}\,F) \boxtimes (\Delta_{\mathcal{I}'}\,\mathrm{colim}\,F') = \Delta_{\mathcal{I}\times\mathcal{I}'}(\mathrm{colim}\,F \otimes \mathrm{colim}\,F')$$

is an isomorphism (see Example 3.7). Let $\iota$ denote the colimiting cocone of $\mathrm{colim}(F \boxtimes F')$. By definition $\Delta\mu^{-1} \circ \iota = \kappa \boxtimes \kappa'$, hence $\iota = \Delta\mu \circ \kappa \boxtimes \kappa'$.

*Naturality of* $\mu^{-1}$. The following square commutes by definition of $\mathrm{colim}(H, \eta)$.

$$\begin{array}{ccc}
(\mathcal{I}, F) & \xrightarrow{(\mathrm{id}, \kappa)} & (\mathcal{I}, \Delta\,\mathrm{colim}\,F) \\
{\scriptstyle (H,\eta)}\downarrow & & \downarrow{\scriptstyle (H, \Delta\,\mathrm{colim}(H,\eta))} \\
(\mathcal{J}, G) & \xrightarrow{(\mathrm{id}, \lambda)} & (\mathcal{J}, \Delta\,\mathrm{colim}\,G)
\end{array}$$

Apply $\boxtimes$ to two such diagrams, obtaining two coinciding cocones $F \boxtimes F' \Rightarrow \Delta_{\mathcal{I}\times\mathcal{I}'}(\mathrm{colim}\,G \otimes \mathrm{colim}\,G')$. By the universal property of $\mathrm{colim}(F \boxtimes F')$, the naturality diagram for $\mu^{-1}$ commutes. ◻

# 4. Colimits of enriched categories

Let $\mathcal{V}$ be cocomplete and let $- \otimes B$ and $A \otimes -$ be cocontinuous. We construct the colimit of a fixed diagram $F\colon \mathcal{I} \to \mathrm{Cat}_{\mathcal{V}}$. Interpret $\mathrm{Ob}\colon \mathrm{Cat}_{\mathcal{V}} \to \mathrm{Set}$ as a forgetful functor. If not explicitly specified, we assume left-associated tensor products by convention.

## 4.1. Strings and their evaluation

Given a diagram $F\colon \mathcal{I} \to \mathrm{Cat}_{\mathcal{V}}$, we will construct the hom-object $(\operatorname{colim} F)(a, b)$ as the colimit of a certain diagram $\mathrm{ev}^F_{a,b}\colon \mathcal{S}_F(a, b) \to \mathcal{V}$. We will introduce the category of strings $\mathcal{S}_F(a, b)$ in Definition 4.5, defining a Cat-enriched category $\mathcal{S}_F$. In Definition 4.8, we define $\mathrm{ev}^F_{a,b}$, lifting $\mathcal{S}_F$ to a $\int[-, \mathcal{V}]$-enriched category $\mathbb{S}_F$. Change of enriching base along the strong monoidal functor colim: $\int[-, \mathcal{V}] \to \mathcal{V}$ from Section 3 yields the $\mathcal{V}$-enriched category $\operatorname{colim}_*(\mathbb{S}_F)$. In Section 4.2, we will prove that $\operatorname{colim}_*(\mathbb{S}_F)$ satisfies the universal property of $\operatorname{colim}(F)$.

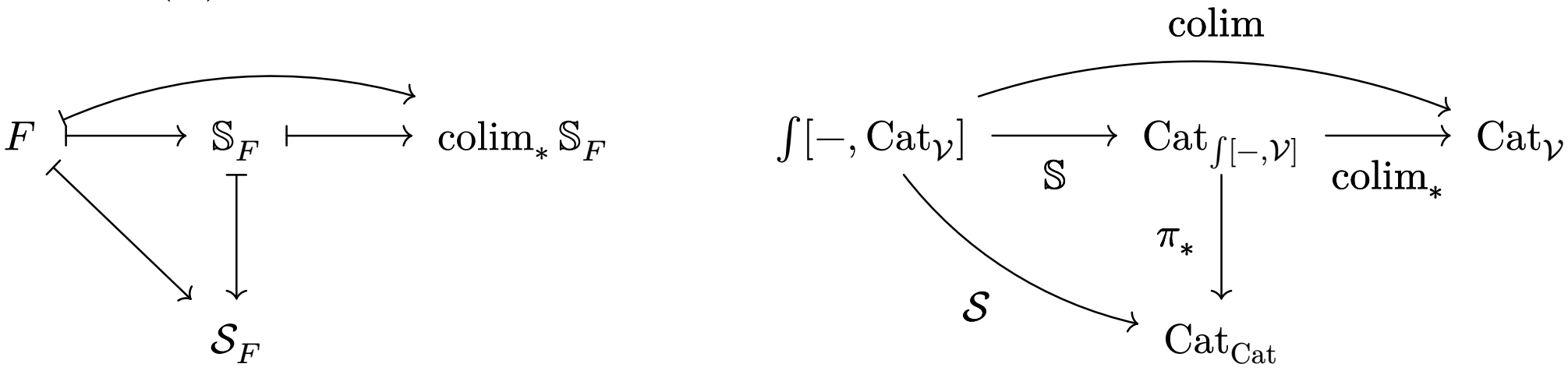

In this section, we will merely prove commutativity of the left diagram. See Appendix B for functoriality in $\int[-, \mathrm{Cat}_{\mathcal{V}}]$.

Definition 4.1 (category of elements) [GR71, §VI.8]. The (covariant) **category of elements** of a functor $G\colon \mathcal{C} \to \mathrm{Set}$ is the category whose

- objects are pairs $(c, e)$ of $c \in \mathcal{C}$ and $e \in G(c)$,
- morphisms $(c, e) \to (c', e')$ are morphisms $f\colon c \to c'$, such that $(Gf)(e) = e'$.

Remark 4.2. The category of elements is a covariant version of the Grothendieck construction $\int F$ from Definition 3.1, considering sets as discrete categories.

Definition 4.3. Let $F\colon \mathcal{I} \to \mathrm{Cat}_{\mathcal{V}}$ be a diagram. The category of elements (see Definition 4.1)

$$\mathcal{E}_F := \int_{i \in \mathcal{I}} \mathrm{Ob}(Fi) \times \mathrm{Ob}(Fi)$$

has as objects ***F*-edges** $(i, x, y)$ with $x, y \in F(i)$ and morphisms $(i, x, y) \to (i', x', y')$ along $f\colon i \to i'$ in $\mathcal{I}$, such that $(Ff)(x) = x'$ and $(Ff)(y) = y'$.

Remark 4.4. Edges constitute a directed graph with vertices $\operatorname{colim}(\mathrm{Ob} \circ F)$, where $(i, x, y)$ is a directed edge from $[x]$ to $[y]$. In the following definition, a "string" is exactly a path in said graph.

Definition 4.5 (category of strings). Let $F\colon \mathcal{I} \to \mathrm{Cat}_{\mathcal{V}}$ be a diagram and let $a, b \in \operatorname{colim}(\mathrm{Ob} \circ F)$. Define the category $\mathcal{S}_F(a, b)$ whose objects are ***F*-strings** $[e_n, ..., e_1]$ of $F$-edges $e_k := (i_k, x_k, y_k) \in \mathcal{E}_F$ from $a = [x_1]$ to $b = [y_n]$ with compatible endpoints $[y_{k-1}] = [x_k] \in \operatorname{colim}(\mathrm{Ob} \circ F)$. If $a = b$, include the **empty string** $\varnothing_a := [\,]$. Denote **concatenation** of strings by $\star_{a,b,c} : \mathcal{S}(b, c) \times \mathcal{S}(a, b) \to \mathcal{S}(a, c)$. Let morphisms in $\mathcal{S}_F(a, b)$ be generated by

- **transports** $[..., f, ...]\colon [..., e, ...] \to [..., e', ...]$, where $f \in \mathcal{E}_F(e, e')$,
- **contractions** $[..., (i, y, z), (i, x, y), ...] \to [..., (i, x, z), ...]$,
- **unit insertions**, such that the corresponding endpoints coincide in $\operatorname{colim}(\mathrm{Ob} \circ F)$. In particular, allow $\varnothing_a \to [(i, x, x)]$ in $\mathcal{S}_F(a, a)$.

$$[..., e_k, ...] \to [..., e_k, (i, x, x), ...] \quad \text{or} \quad [..., e_k, ...] \to [..., (j, y, y), e_k, ...]$$

$$\begin{array}{ccc} t \star s & \xrightarrow{\mathrm{id} \star \alpha} & t \star s' \\ {\scriptstyle \beta \star \mathrm{id}} \downarrow & & \downarrow {\scriptstyle \beta \star \mathrm{id}} \\ t' \star s & \xrightarrow[\mathrm{id} \star \alpha]{} & t' \star s' \end{array}$$

To ensure functoriality of string concatenation, quotient by the **interchange** relation on the right for any generators $\alpha\colon s \to s'$ and $\beta\colon t \to t'$.

Here, we encode unit maps in $\mathcal{S}_F(a, b)$ directly, so that our construction of $\operatorname{colim}(F)$ in Theorem 4.14 arises via a change of enriching base. In Definition 4.17, we leave out empty strings and unit insertions and define unit maps $\mathbb{1} \to (\operatorname{colim} F)(a, a)$ in Corollary 4.21 externally.

The generating morphisms in $\mathcal{S}_F(a, b)$ are motivated by the universal cocone of $\operatorname{colim}(F)$, later defined in Lemma 4.15. Transports ensure that the relevant triangles in said cocone commute. Contractions encode composition coherence of the cocones' components, while empty strings and unit insertions $\varnothing_a \to [(i,x,x)]$ encode unit coherence. The interchange relation makes independent local modifications of a string commute.

Remark 4.6. We can organize $\mathcal{S}_F$ into a Cat-enriched category with objects $\operatorname{colim}(\mathrm{Ob}\circ F)$ and hom-categories $\mathcal{S}(a,b)$. Composition is given by the concatenation bifunctors $\star_{a,b,c}$ and unit maps $* \to \mathcal{S}(a,a)$ by $* \mapsto \varnothing_a$. It is easily checked that $\mathcal{S}_F$ satisfies associativity and unitality.

Remark 4.7. Alternatively, one can extend the graph from Remark 4.4 to a *computad* [Str76, §2] by adjoining generating 2-cells corresponding to transports, contractions and unit insertions respectively. This computad freely generates the strict 2-category $\mathcal{S}_F$ from Remark 4.6. In particular, the interchange relation is automatic.

$$a \underset{(i',x',y')}{\overset{(i,x,y)}{\rightrightarrows}} b \;(\Downarrow) \qquad a \xrightarrow{(i,x,y)} b \xrightarrow{(i,y,z)} c \;\Downarrow\; a \xrightarrow{(i,x,z)} c \qquad a \;\; \varnothing_a \Rightarrow (i,x,x)$$

Definition 4.8. Define the ***F*-string evaluation** functor $\mathrm{ev}^F_{a,b}\colon \mathcal{S}_F(a,b) \to \mathcal{V}$ on non-empty strings by

$$[(i_n,x_n,y_n),\ldots,(i_1,x_1,y_1)] \mapsto F(i_n)(x_n,y_n)\otimes\cdots\otimes F(i_1)(x_1,y_1)$$

and $\varnothing_a \mapsto \mathbb{1}$ if $a=b$. Send edgewise transports to $Ff_n\otimes\cdots\otimes Ff_1$, contractions to $\circ\colon F(i)(y,z)\otimes F(i)(x,y) \to F(i)(x,z)$ and unit insertions via the unitors' inverses to

$$\cdots\otimes\cdots \cong \cdots\otimes\mathbb{1}\otimes\cdots \xrightarrow{\mathrm{id}\otimes j_x\otimes\mathrm{id}} \cdots\otimes F(i)(x,x)\otimes\cdots$$

For generating morphisms $\alpha\colon s\to s'$ in $\mathcal{S}_F(a,b)$ and $\beta\colon t\to t'$ in $\mathcal{S}_F(b,c)$, the two sides of their interchange square evaluate, up to canonical reassociation and unitor isomorphisms, to

$$(\mathrm{ev}(\beta)\otimes\mathrm{id})\circ(\mathrm{id}\otimes\mathrm{ev}(\alpha)) = \mathrm{ev}(\beta)\otimes\mathrm{ev}(\alpha) = (\mathrm{id}\otimes\mathrm{ev}(\alpha))\circ(\mathrm{ev}(\beta)\otimes\mathrm{id}),$$

by functoriality of the tensor product $\otimes$ and Mac Lane coherence [Mac98, §VII.2]. Hence $\mathrm{ev}^F_{a,c}$ is well-defined on the interchange quotient.

Remark 4.9. We could have added additional relations in Definition 4.5, corresponding to functoriality of transport, as well as monoidal and enriched coherence axioms. Equivalently, we could have adjoined 3-cells to the computad from Remark 4.7 and quotiented the free 2-category by the congruence generated by the boundary pairs of said 3-cells.

Lemma 4.10. The associator and the unitors of $\mathcal{V}$ define a family of natural isomorphisms

$$\theta_{a,b,c}\colon \mathrm{ev}^F_{b,c}\boxtimes\mathrm{ev}^F_{a,b} \overset{\cong}{\Longrightarrow} \mathrm{ev}^F_{a,c}\big(-\star_{a,b,c}-\big).$$

Proof. Write $\mathrm{ev} := \mathrm{ev}^F$ and suppress the superscript endpoints. For non-empty strings $s\in\mathcal{S}_F(a,b)$ and $s'\in\mathcal{S}_F(b,c)$, iterative application of the associator yields $\mathrm{ev}(s')\otimes\mathrm{ev}(s)\cong\mathrm{ev}(s'\star s)$. Applying the left unitor gives $\mathbb{1}\otimes\mathrm{ev}(s)\cong\mathrm{ev}(s)$ and the right unitor gives $\mathrm{ev}(s')\otimes\mathbb{1}\cong\mathrm{ev}(s')$.

*Naturality.* Each isomorphism $\theta_{a,b,c}$ is natural in transports by naturality of the associator, natural in contractions by associativity of enriched composition and Mac Lane coherence, and natural in unit insertions by enriched unitality, naturality of the associator and inverse unitors and Mac Lane coherence. □

Definition 4.11. In the sense of Remark 4.12, lift $\mathcal{S}_F$ to a category $\mathbb{S}_F$ enriched in $\int[-,\mathcal{V}]$ with the monoidal structure from Lemma 3.5, with

- objects $\mathrm{Ob}(\mathbb{S}_F) := \mathrm{Ob}(\mathcal{S}_F) = \operatorname{colim}(\mathrm{Ob}\circ F)$,
- hom-diagrams $\mathbb{S}_F(a,b) := \big(\mathcal{S}_F(a,b),\mathrm{ev}^F_{a,b}\big)$,
- composition $\circ_{a,b,c} := \big(\star_{a,b,c},\theta_{a,b,c}\big)$ via string concatenation $\star_{a,b,c}$ and $\theta_{a,b,c}$ from Lemma 4.10,
- and unit maps $\big(j^{\mathcal{S}}_a,\mathrm{id}\big)\colon (*\mapsto\mathbb{1}) \to \mathbb{S}_F(a,a)$.

In Lemma A.3 of Appendix A, we verify the coherence axioms for the enriched category $\mathbb{S}_F$ by untangling compositions in $\int[-,\mathcal{V}]$ and using the enriched axioms of $\mathcal{S}_F$.

Remark 4.12. Lemma 2.5 applied to the strict monoidal functor $\pi\colon \int[-,\mathcal{V}] \to \mathrm{Cat}$ sending $(\mathcal{I}, F) \mapsto \mathcal{I}$ induces a functor $\pi_*\colon \mathrm{Cat}_{\int[-,\mathcal{V}]} \to \mathrm{Cat}_{\mathrm{Cat}}$ and $\pi_*(\mathbb{S}_F)$ is the Cat-enriched $\mathcal{S}_F$ from Remark 4.6.

Description 4.13. Change of enriching base along the strong monoidal functor colim: $\int[-,\mathcal{V}] \to \mathcal{V}$ from Proposition 3.8 induces a functor $\mathrm{colim}_*\colon \mathrm{Cat}_{\int[-,\mathcal{V}]} \to \mathrm{Cat}_{\mathcal{V}}$ by Lemma 2.5. Unwinding $\mathrm{colim}_*(\mathbb{S}_F)$ we have objects $\mathrm{colim}(\mathrm{Ob} \circ F)$ and hom-objects $\mathrm{colim}_{\mathcal{S}_F(a,b)}\left(\mathrm{ev}^F_{a,b}\right)$, with the following composition and unit maps.

$$\mathrm{colim}\, \mathbb{S}_F(b,c) \otimes \mathrm{colim}\, \mathbb{S}_F(a,b) \overset{\mu}{\cong} \mathrm{colim}(\mathbb{S}_F(b,c) \boxtimes \mathbb{S}_F(a,b)) \xrightarrow{\mathrm{colim}(\star,\theta)} \mathrm{colim}\, \mathbb{S}_F(a,c)$$

$$\mathbb{1} \overset{\eta}{\cong} \mathrm{colim}(* \mapsto \mathbb{1}) \xrightarrow{\mathrm{colim}(* \mapsto \varnothing,\, \mathrm{id})} \mathrm{colim}\, \mathbb{S}_F(a,a)$$

## 4.2. The colimit theorem

Theorem 4.14. Let $\mathcal{V}$ be cocomplete and let $\otimes$ be cocontinuous in each argument. Then $\mathrm{colim}_*(\mathbb{S}_F)$ is a colimit of the small diagram $F\colon \mathcal{I} \to \mathrm{Cat}_{\mathcal{V}}$.

Proof. In Lemma 4.15, we construct a cocone $\iota\colon F \Rightarrow \Delta\, \mathrm{colim}_*\, \mathbb{S}_F$ by embedding $F(i) \to \mathrm{colim}_*(\mathbb{S}_F)$ via singleton strings $[(i,x,y)]$. In Lemma 4.16, we show universality of $\iota$. Given a cocone $\eta\colon F \Rightarrow \Delta\mathcal{D}$, we construct an enriched functor $\tilde{\eta}\colon \mathrm{colim}_*(\mathbb{S}_F) \to \mathcal{D}$ on objects via the universal property of $\mathrm{colim}(\mathrm{Ob} \circ F)$ and on hom-objects via the universal property of $\mathrm{colim}\, \mathbb{S}_F(a,b)$. We obtain a cocone $\mathbb{S}_F(a,b) \Rightarrow \Delta\mathcal{D}(\tilde{\eta}a, \tilde{\eta}b)$ by applying $\eta$ to each edge, followed by composition in $\mathcal{D}$. Uniqueness of $\tilde{\eta}$ follows, since any string decomposes into a concatenation of singleton strings. □

Lemma 4.15. Let $\mathrm{colim}_*(\mathbb{S}_F)$ be as in Theorem 4.14. The following components $\iota_i\colon F(i) \to \mathrm{colim}_*(\mathbb{S}_F)$ define a cocone $\iota\colon F \Rightarrow \Delta\, \mathrm{colim}_*(\mathbb{S}_F)$.

- on objects $\iota_i(x) := [x] \in \mathrm{colim}(\mathrm{Ob} \circ F)$,
- on hom-objects $(\iota_i)_{x,y}\colon F(i)(x,y) \to \mathrm{colim}\, \mathbb{S}_F([x],[y])$ by sending $[(i,x,y)] \in \mathcal{S}([x],[y])$ via the colimiting cocone.

Proof. Write $\mathbb{S} := \mathbb{S}_F$. First show that the $\iota_i$ are enriched functors.

*Composition.* The diagonal arrow is a component of the cocone of $\mathrm{colim}\, \mathbb{S}([x],[z])$, thus the upper right triangle commutes. The lower left triangle commutes, since $\mathrm{colim}(\mathbb{S}([y],[z]) \boxtimes \mathbb{S}([x],[y])) \to \mathrm{colim}\, \mathbb{S}([x],[z])$ was defined as the unique morphism between the cocone of $\mathrm{colim}(\mathbb{S}([y],[z]) \boxtimes \mathbb{S}([x],[y]))$ and the cocone $\mathbb{S}([y],[z]) \boxtimes \mathbb{S}([x],[y]) \Rightarrow \Delta\, \mathrm{colim}\, \mathbb{S}([x],[z])$.

$$\begin{array}{ccc} F(i)(y,z) \otimes F(i)(x,y) & \xrightarrow{\circ} & F(i)(x,z) \\ \downarrow{\scriptstyle (\iota_i)_{y,z} \otimes (\iota_i)_{x,y}} & \searrow & \downarrow{\scriptstyle (\iota_i)_{x,z}} \\ \mathrm{colim}\, \mathbb{S}([y],[z]) \otimes \mathrm{colim}\, \mathbb{S}([x],[y]) & \xrightarrow[\circ]{} & \mathrm{colim}\, \mathbb{S}([x],[z]) \end{array}$$

*Units.* The lower left triangle is part of the cocone of $\mathrm{colim}\, \mathbb{S}([x],[x])$. The upper right triangle commutes by the definition of $\mathrm{colim}(* \mapsto \varnothing, \mathrm{id})$.

$$\begin{array}{ccc} \mathbb{1} & \xrightarrow{\cong} & \mathrm{colim}(* \mapsto \mathbb{1}) \\ \downarrow{\scriptstyle j_x} & \searrow & \downarrow{\scriptstyle \mathrm{colim}(* \mapsto \varnothing,\, \mathrm{id})} \\ F(i)(x,x) & \xrightarrow[(\iota_i)_{x,x}]{} & \mathrm{colim}\, \mathbb{S}([x],[x]) \end{array}$$

*Cocone.* The $\iota_i$ form a cocone using the cocones of $\operatorname{colim} \mathbb{S}(a,b)$ over string morphisms $[(i,x,y)] \mapsto [(j,(Ff)(x),(Ff)(y))]$.

$$\begin{array}{ccc} F(i)(x,y) & \xrightarrow{F(f)_{x,y}} & F(j)((Ff)(x),(Ff)(y)) \\ & {\scriptstyle (\iota_i)_{x,y}} \searrow \quad \swarrow {\scriptstyle (\iota_j)_{(Ff)(x),(Ff)(y)}} & \\ & \operatorname{colim} \mathbb{S}([x],[y]) & \end{array}$$

□

Lemma 4.16. The cocone $\iota\colon F \Rightarrow \Delta \operatorname{colim}_*(\mathbb{S}_F)$ from Lemma 4.15 is universal.

Proof. Write $\mathbb{S} := \mathbb{S}_F$. Let $\eta\colon F \Rightarrow \Delta\mathcal{D}$ and find a unique enriched functor $\tilde{\eta}\colon \operatorname{colim}_* \mathbb{S} \to \mathcal{D}$ with $\eta = \Delta\tilde{\eta} \circ \iota$.

*Objects.* By the universal property of $\operatorname{colim}(\mathrm{Ob} \circ F)$, the cocone $\mathrm{Ob}\,\eta\colon \mathrm{Ob} \circ F \Rightarrow \Delta\,\mathrm{Ob}(\mathcal{D})$ induces a unique function $\mathrm{Ob}(\operatorname{colim}_* \mathbb{S}) \to \mathrm{Ob}(\mathcal{D})$.

*Hom-objects.* By constructing a cocone $\mathbb{S}(a,b) \Rightarrow \Delta\mathcal{D}(\tilde{\eta}a, \tilde{\eta}b)$, the universal property of $\operatorname{colim} \mathbb{S}(a,b)$ induces a unique morphism $\tilde{\eta}_{a,b}$. For non-empty strings, define components

$$F(j)(x_n,y) \otimes \cdots \otimes F(i)(x,y_1) \xrightarrow{\eta_j \otimes \cdots \otimes \eta_i} \mathcal{D}\big(\eta_j x_n, \eta_j y\big) \otimes \cdots \otimes \mathcal{D}(\eta_i x, \eta_i y_1) \xrightarrow{\circ^{n-1}} \mathcal{D}\big(\eta_i x, \eta_j y\big). \tag{2}$$

In $\mathbb{S}(a,a) \Rightarrow \Delta\mathcal{D}(\tilde{\eta}a, \tilde{\eta}a)$, let the component $\mathbb{1} \to \mathcal{D}(\tilde{\eta}a, \tilde{\eta}a)$ be the unit map. Edgewise transports $Ff_n \otimes \cdots \otimes Ff_1$ in $\mathbb{S}(a,b)$ are respected, since $\eta$ is a cocone and $\otimes$ is a functor. Contractions $F(i)(y,z) \otimes F(i)(x,y) \xrightarrow{\circ} F(i)(x,z)$ are respected by composition coherence of the enriched functor $\eta_i$. Unit insertions $j_x\colon \mathbb{1} \to F(i)(x,x)$ are respected by unit coherence of $\eta_i$.

*Units.* Since the component $\mathbb{1} \to \mathcal{D}(\tilde{\eta}a, \tilde{\eta}a)$ in the cocone $\mathbb{S}(a,a) \Rightarrow \Delta\mathcal{D}(\tilde{\eta}a, \tilde{\eta}a)$ was defined as the unit map, the unital diagram commutes by the universal property of $\operatorname{colim} \mathbb{S}(a,a)$.

*Composition.* Let $s \in \mathcal{S}_F(a,b)$. Recall from Definition 4.8 that $\mathrm{ev}^F(\varnothing) = \mathbb{1}$. The top square commutes by definition of $\theta$ and naturality of the right unitor. The bottom triangle commutes by the unital diagram of $\mathcal{D}$. Composition with the empty string is analogous.

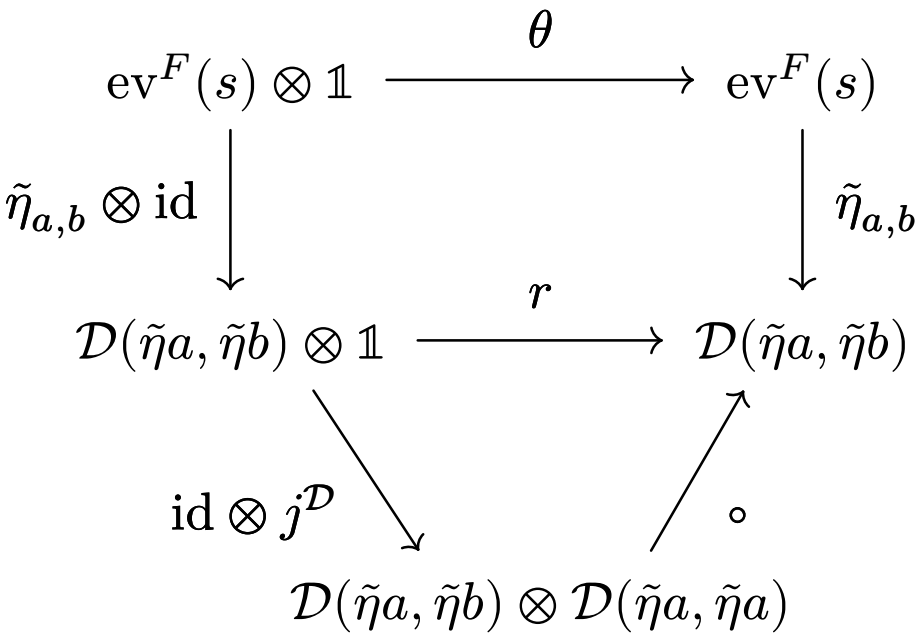

For non-empty strings, the top square commutes by naturality of the associator $\alpha$. The bottom square commutes by the associativity diagram of $\mathcal{D}$.

$$\begin{array}{ccc} \Big(\bigotimes_{k=m+n}^{m+1} F(i_k)(x_k,y_k)\Big) \otimes \Big(\bigotimes_{k=m}^{1} F(i_k)(x_k,y_k)\Big) & \xrightarrow{\theta} & \bigotimes_{k=m+n}^{1} F(i_k)(x_k,y_k) \\ \Big\downarrow {\scriptstyle \left(\bigotimes \eta_{i_k}\right) \otimes \left(\bigotimes \eta_{i_k}\right)} & & \Big\downarrow {\scriptstyle \bigotimes \eta_{i_k}} \\ \Big(\bigotimes_{k=m+n}^{m+1} \mathcal{D}(\eta x_k, \eta y_k)\Big) \otimes \Big(\bigotimes_{k=m}^{1} \mathcal{D}(\eta x_k, \eta y_k)\Big) & \longrightarrow & \bigotimes_{k=m+n}^{1} \mathcal{D}(\eta x_k, \eta y_k) \\ \Big\downarrow {\scriptstyle \circ^{n-1} \otimes \circ^{m-1}} & & \Big\downarrow {\scriptstyle \circ^{m+n-1}} \\ \mathcal{D}(\eta x_{m+1}, \eta y_{m+n}) \otimes \mathcal{D}(\eta x_1, \eta y_m) & \xrightarrow{\circ} & \mathcal{D}(\eta x_1, \eta y_{m+n}) \end{array}$$

Let $a \coloneqq [x_1]$, $b \coloneqq [y_m] = [x_{m+1}]$ and $c \coloneqq [y_{m+n}]$. The upper right and lower left compositions define the same cocone $\mathbb{S}(b,c) \boxtimes \mathbb{S}(a,b) \Rightarrow \Delta\mathcal{D}(\tilde\eta a, \tilde\eta c)$. By the universal property, the induced morphisms agree.

$$
\begin{array}{ccccc}
\operatorname{colim}\mathbb{S}(b,c) \otimes \operatorname{colim}\mathbb{S}(a,b) & \xrightarrow{\cong} & \operatorname{colim}(\mathbb{S}(b,c) \boxtimes \mathbb{S}(a,b)) & \xrightarrow{\operatorname{colim}(\star,\theta)} & \operatorname{colim}\mathbb{S}(a,c) \\
 & \searrow^{\tilde\eta_{b,c} \otimes \tilde\eta_{a,b}} & \downarrow & & \downarrow{\tilde\eta_{a,c}} \\
 & & \mathcal{D}(\tilde\eta b, \tilde\eta c) \otimes \mathcal{D}(\tilde\eta a, \tilde\eta b) & \xrightarrow{\circ} & \mathcal{D}(\tilde\eta a, \tilde\eta c)
\end{array}
$$

The left triangle commutes, since by definition of $\operatorname{colim} F \otimes \operatorname{colim} G \cong \operatorname{colim}(F \boxtimes G)$ from Proposition 3.8, the following diagram in $[\mathcal{S}(b,c) \times \mathcal{S}(a,b), \mathcal{V}]$ commutes.

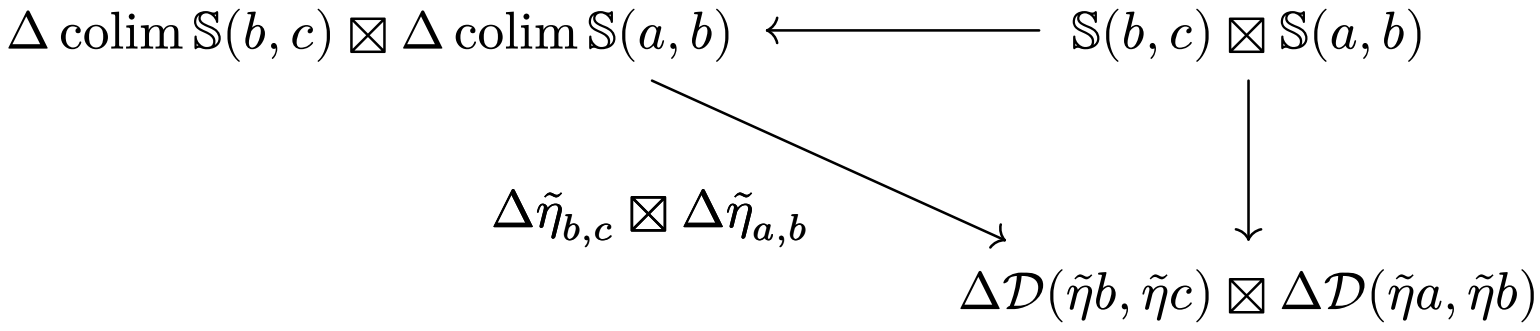

*Uniqueness.* By the universal property of $\operatorname{colim}(\mathrm{Ob} \circ F)$, any enriched functor $\varphi\colon \operatorname{colim}_* \mathbb{S} \to \mathcal{D}$ between cocones $\iota\colon F \Rightarrow \Delta \operatorname{colim}_* \mathbb{S}$ and $\eta\colon F \Rightarrow \Delta\mathcal{D}$ must agree with $\tilde\eta$ on objects. By the universal property of $\operatorname{colim}\mathbb{S}(a,b)$, the morphism $\varphi^{a,b}$ between hom-objects is uniquely determined by a cocone $\mathbb{S}(a,b) \Rightarrow \Delta\mathcal{D}(\tilde\eta a, \tilde\eta b)$.

Let $\kappa^{a,b}$ denote the cocone of $\operatorname{colim}\mathbb{S}(a,b)$. For non-empty strings, the cocone's component must match Equation 2, since any non-empty string $[e_n, \ldots, e_1]$ can be decomposed into a concatenation $[e_n] \star \cdots \star [e_1]$ and the enriched functor $\varphi\colon \operatorname{colim}_* \mathbb{S} \to \mathcal{D}$ must respect composition. For a length-one string $s \coloneqq [(i,x,y)]$, remember $\kappa_s = (\iota_i)_{x,y}$ from Lemma 4.15, thus $\varphi_{[x],[y]} \circ \kappa_s = (\varphi \circ \iota_i)_{x,y} = (\eta_i)_{x,y}$.

Inductively, each component $\varphi_{a,c} \circ \kappa^{a,c}$ is uniquely determined for a non-empty string in $\mathcal{S}_F(a,c)$, which can be decomposed into a concatenation of non-empty strings $s \in \mathcal{S}_F(a,b)$ and $s' \in \mathcal{S}_F(b,c)$.

$$
\begin{array}{ccc}
\mathrm{ev}^F_{b,c}(s') \otimes \mathrm{ev}^F_{a,b}(s) & \xleftarrow{\theta^{-1}} & \mathrm{ev}^F_{a,c}(s' \star s) \\
\downarrow{\kappa^{b,c} \otimes \kappa^{a,b}} & & \downarrow{\kappa^{a,c}} \\
\operatorname{colim}\mathbb{S}(b,c) \otimes \operatorname{colim}\mathbb{S}(a,b) & \xrightarrow{\circ} & \operatorname{colim}\mathbb{S}(a,c) \\
\downarrow{\varphi_{b,c} \otimes \varphi_{a,b}} & & \downarrow{\varphi_{a,c}} \\
\mathcal{D}(\tilde\eta b, \tilde\eta c) \otimes \mathcal{D}(\tilde\eta a, \tilde\eta b) & \xrightarrow{\circ} & \mathcal{D}(\tilde\eta a, \tilde\eta c)
\end{array}
$$

The upper triangle is both the unit diagram of $\iota_i$ from Lemma 4.15 and the universal cocone of $\operatorname{colim}\mathbb{S}(a,a)$ for $\varnothing_a \to [(i,x,x)]$. The component $\mathbb{1} \to \mathcal{D}(\tilde\eta a, \tilde\eta a)$ of the empty string $\varnothing_a$ is forced to be the unit map, since the enriched functor $F(i) \to \mathcal{D}$ respects units.

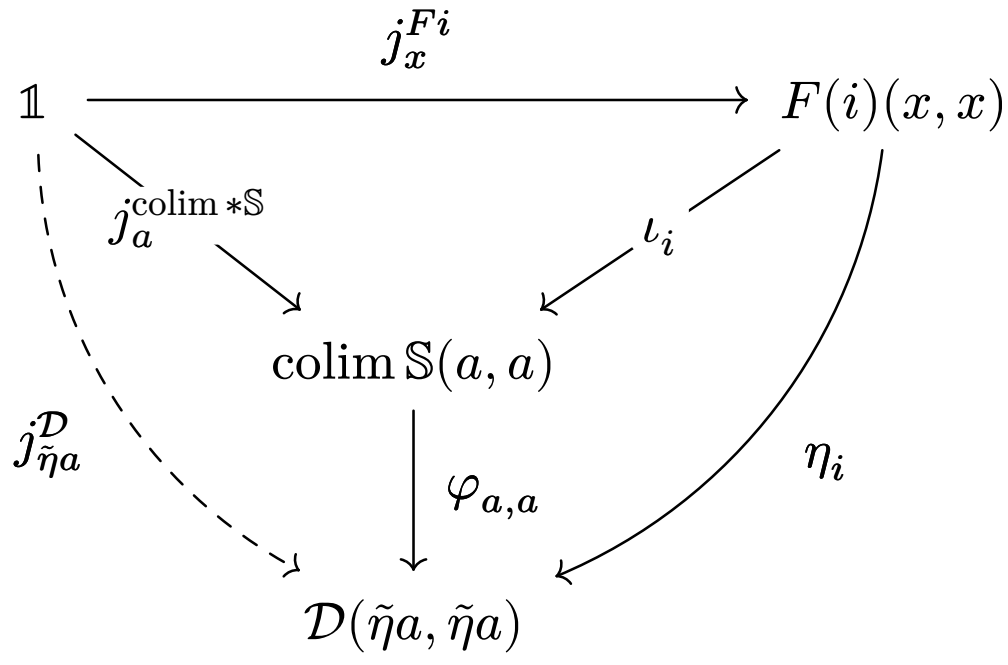

□

## 4.3. Category of non-empty strings

We simplify $\mathcal{S}_F$ from Definition 4.5 to obtain a subcategory of non-empty strings $\overline{\mathcal{S}}_F$. This requires us to define unit maps $\mathbb{1} \to (\operatorname{colim} F)(a,a)$ externally in Corollary 4.21.

Definition 4.17 (category of non-empty strings). Let $F\colon \mathcal{I} \to \mathrm{Cat}_{\mathcal{V}}$ be a diagram and let $a, b \in \operatorname{colim}(\mathrm{Ob} \circ F)$. Define the subcategory $\overline{\mathcal{S}}_F(a,b) \subset \mathcal{S}_F(a,b)$ of non-empty strings, with morphisms generated by transports $[..., f, ...]$ and contractions $[..., (i,y,z), (i,x,y), ...] \to [..., (i,x,z), ...]$, modulo the interchange relation from Definition 4.5.

Lemma 4.18. Let $\mathcal{V}$ be locally small, and let $G\colon \mathcal{J} \to \mathcal{V}$ and $H\colon \mathcal{K} \to \mathcal{V}$ be small diagrams. A natural isomorphism $\mathrm{Cocone}(G,-) \simeq \mathrm{Cocone}(H,-)$ induces an isomorphism $\operatorname{colim}_{\mathcal{J}}(G) \cong \operatorname{colim}_{\mathcal{K}}(H)$.

Proof. Let $A \coloneqq \operatorname{colim}(G)$ and $B \coloneqq \operatorname{colim}(H)$. Let $\alpha$ denote the composite of natural isomorphisms

$$\mathcal{V}(A,-) \simeq \mathrm{Cocone}(G,-) \simeq \mathrm{Cocone}(H,-) \simeq \mathcal{V}(B,-).$$

By [Rie16, Proposition 2.3.1], $\alpha$ determines an isomorphism $A \cong B$. Explicitly, the Yoneda lemma gives a bijection

$$\mathrm{Nat}(\mathcal{V}(A,-), \mathcal{V}(B,-)) \cong \mathcal{V}(B,A),$$

where $\alpha$ corresponds to $g \coloneqq \alpha_A(\mathrm{id}_A)\colon B \to A$ and $\alpha^{-1}$ corresponds to $f \coloneqq \alpha_B^{-1}(\mathrm{id}_B)\colon A \to B$. Insert $\mathrm{id}_A$ in the left diagram to obtain the right diagram, where $f \circ g = \mathrm{id}_B$. The proof that $g \circ f = \mathrm{id}_A$ is analogous.

$$\begin{array}{ccc} \mathcal{V}(A,A) & \xrightarrow{\alpha_A} & \mathcal{V}(B,A) \\ {\scriptstyle f\circ -}\downarrow & & \downarrow{\scriptstyle f\circ -} \\ \mathcal{V}(A,B) & \xrightarrow{\alpha_B} & \mathcal{V}(B,B) \end{array} \qquad \begin{array}{ccc} \mathrm{id}_A & \overset{\alpha_A}{\longmapsto} & g \\ {\scriptstyle f\circ -}\downarrow & & \downarrow{\scriptstyle f\circ -} \\ f & \overset{\alpha_B}{\longmapsto} & \mathrm{id}_B = f \circ g \end{array}$$

□

Lemma 4.19. Let $F\colon \mathcal{I} \to \mathrm{Cat}_{\mathcal{V}}$ be a diagram and let $\kappa\colon \overline{\mathrm{ev}}_F^{a,a} \Rightarrow \Delta v$ be a cocone. Then the composite

$$\mathbb{1} \xrightarrow{j_x} F(i)(x,x) \xrightarrow{\kappa} v$$

is independent of the representative of $[x] \in \operatorname{colim}(\mathrm{Ob} \circ F)$.

Proof. Any edge $(j,y,y) \in \mathcal{E}_F$ as in Definition 4.3 with $[y] = [x] \in \operatorname{colim}(\mathrm{Ob} \circ F)$ is connected to the edge $(i,x,x)$ by an underlying zigzag $\mathcal{Z} \to \mathcal{I}$ from $i$ to $j$. Since enriched functors respect units, we get a cone $\Delta\mathbb{1} \Rightarrow (\mathcal{Z} \to \mathcal{S}_F(a,a) \to \mathcal{V})$ and all components of the composite $\Delta\mathbb{1} \Rightarrow \Delta v$ coincide.

$$\begin{array}{ccc} \mathbb{1} & \xrightarrow{j_x} & F(i)(x,x) \\ {\scriptstyle j_y}\downarrow & \rightsquigarrow & \downarrow{\scriptstyle \kappa} \\ F(j)(y,y) & \xrightarrow[\kappa]{} & v \end{array}$$

□

Lemma 4.20. Let $\mathcal{V}$ be locally small and let $F\colon \mathcal{I} \to \mathrm{Cat}_{\mathcal{V}}$ be a small diagram. The inclusion $\iota\colon \overline{\mathcal{S}}_F(a,b) \hookrightarrow \mathcal{S}_F(a,b)$ induces an isomorphism $\operatorname{colim}\big(\mathrm{ev}_{a,b}^F\big) \cong \operatorname{colim}\big(\mathrm{ev}_{a,b}^F \circ \iota\big)$.

Proof. Write $\mathrm{ev} \coloneqq \mathrm{ev}_{a,b}^F$. By Lemma 4.18, it suffices to construct a natural isomorphism in $v \in \mathcal{V}$ between cocone categories

$$\mathrm{Cocone}(\mathrm{ev}, v) \underset{\rho_v}{\overset{\varepsilon_v}{\rightleftarrows}} \mathrm{Cocone}(\mathrm{ev} \circ \iota, v),$$

Restrict a cocone $\lambda\colon \mathrm{ev} \Rightarrow \Delta v$ along $\iota$ by whiskering $\lambda\iota \eqqcolon \rho_v(\lambda)$. This restriction is clearly natural in $v$.

$$\mathrm{ev} \circ \iota \overset{\lambda\iota}{\Longrightarrow} (\Delta v) \circ \iota = \Delta v$$

We now define the inverse by extending a cocone $\kappa$ under $\mathrm{ev} \circ \iota$ with tip $v$ to ev. For non-empty strings $s \in \overline{\mathcal{S}}_F(a,b)$, let $\varepsilon_v(\kappa)_s := \kappa_s$.

*Empty string.* If $a = b$, let $\varepsilon_v(\kappa)_{\varnothing_a}$ be some composite $\mathbb{1} \xrightarrow{j_x} F(i)(x,x) \xrightarrow{\kappa} v$ for $[x] = a \in \mathrm{colim}(\mathrm{Ob} \circ F)$. This composite is independent of the representative $x$ by Lemma 4.19.

*Unit insertion.* For unit insertions $u\colon [..., (i,x,y), ...] \to [..., (i,x,y), (i,x,x), ...]$, where the inserted edge comes from the same $i \in \mathcal{I}$, the following triangle commutes by the right unital diagram and the definition of $\mathrm{ev}(u)$ (see Definition 4.8). Hence, the cocone is natural in unit insertions.

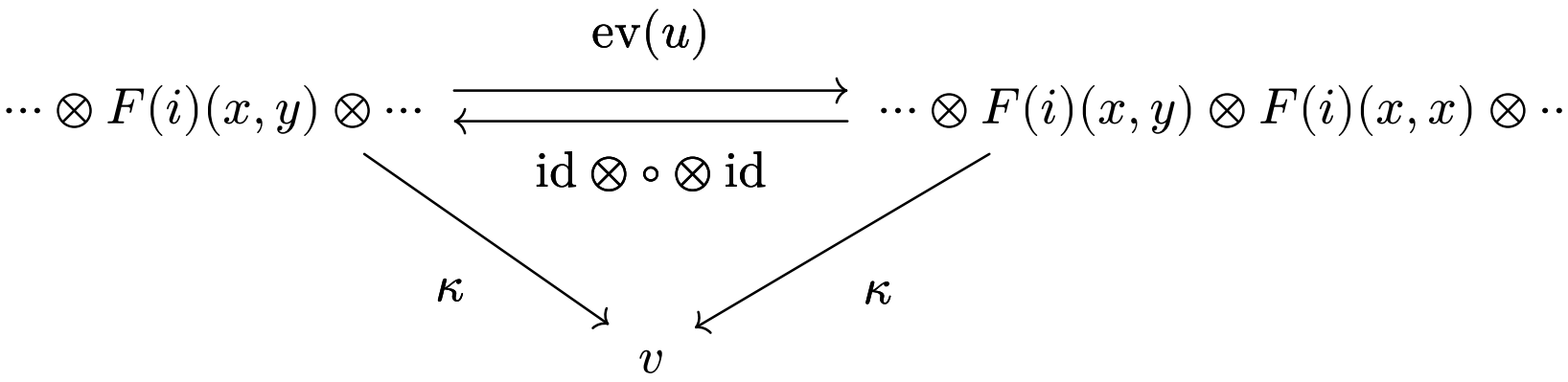

For unit insertions $[..., (i,x,y), ...] \to [..., (i,x,y), (j,x',x'), ...]$ where $i \neq j$, first connect $F(j)(x',x')$ to $F(i)(x,x)$ via a zigzag as in the proof of Lemma 4.19, then proceed as in the above case of $i = j$.

*Inverses.* If $a = b$, extending the restriction of some $\lambda\colon \mathrm{ev} \Rightarrow \Delta v$ recovers $\lambda_{\varnothing_a} = \varepsilon_v(\rho_v(\lambda))_{\varnothing_a}$ by naturality of $\lambda$ in $\varnothing_a \to [(i,x,x)]$. □

Corollary 4.21. Let $\mathcal{V}$ be cocomplete and locally small, let $\otimes$ be cocontinuous in each argument, and let $F\colon \mathcal{I} \to \mathrm{Cat}_{\mathcal{V}}$ be a small diagram. A colimit of $F$ is given by the $\mathcal{V}$-enriched category, where

- objects are $\mathrm{colim}(\mathrm{Ob} \circ F)$,
- for composites $\overline{\mathrm{ev}}_F^{a,b}\colon \overline{\mathcal{S}}_F(a,b) \hookrightarrow \mathcal{S}_F(a,b) \xrightarrow{\mathrm{ev}} \mathcal{V}$, hom-objects are $\mathrm{colim}\big(\overline{\mathrm{ev}}_F^{a,b}\big)$,
- a unit map $\mathbb{1} \to F(i)(x,x) \to \mathrm{colim}(\overline{\mathrm{ev}}_F^{a,a})$ is the composite of the unit map $j_x$ of $F(i)$ and the colimiting cocone.

Let $\overline{\star}_{a,b,c} : \overline{\mathcal{S}}_F(b,c) \times \overline{\mathcal{S}}_F(a,b) \to \overline{\mathcal{S}}_F(a,c)$ denote the restriction of the string concatenation $\star_{a,b,c}$ from Remark 4.6 and let $\overline{\theta}_{a,b,c}$ denote the restriction of $\theta_{a,b,c}$ from Lemma 4.10. Via the natural isomorphism $\mu$ from Proposition 3.8, define composition $\circ_{a,b,c}$ by

$$\mathrm{colim}\big(\overline{\mathrm{ev}}_F^{b,c}\big) \otimes \mathrm{colim}\big(\overline{\mathrm{ev}}_F^{a,b}\big) \overset{\mu}{\cong} \mathrm{colim}\big(\overline{\mathrm{ev}}_F^{b,c} \boxtimes \overline{\mathrm{ev}}_F^{a,b}\big) \xrightarrow{\mathrm{colim}(\overline{\star},\overline{\theta})} \mathrm{colim}(\overline{\mathrm{ev}}_F^{a,c}).$$

Proof. We derive this corollary from Theorem 4.14, Lemma 4.19 and Lemma 4.20. Denote the inclusion $\iota_{a,b}\colon \overline{\mathcal{S}}_F(a,b) \hookrightarrow \mathcal{S}_F(a,b)$.

*Composition.* The left square commutes by naturality of $\mu$ from Proposition 3.8. The right square commutes by functoriality of colim on $\int[-,\mathcal{V}]$ (see Proposition 3.3) and $\star_{a,b,c} \circ \big(\iota_{b,c} \times \iota_{a,b}\big) = \iota_{a,c} \circ \overline{\star}_{a,b,c}$

$$\begin{array}{ccccc}
\mathrm{colim}\big(\overline{\mathrm{ev}}_F^{b,c}\big) \otimes \mathrm{colim}\big(\overline{\mathrm{ev}}_F^{a,b}\big) & \xrightarrow{\mu} & \mathrm{colim}\big(\overline{\mathrm{ev}}_F^{b,c} \boxtimes \overline{\mathrm{ev}}_F^{a,b}\big) & \xrightarrow{\mathrm{colim}(\overline{\star},\overline{\theta})} & \mathrm{colim}(\overline{\mathrm{ev}}_F^{a,c}) \\
\downarrow{\scriptstyle \mathrm{colim}(\iota,\mathrm{id}) \otimes \mathrm{colim}(\iota,\mathrm{id})} & & \downarrow{\scriptstyle \mathrm{colim}(\iota \times \iota,\mathrm{id})} & & \downarrow{\scriptstyle \mathrm{colim}(\iota,\mathrm{id})} \\
\mathrm{colim}\big(\mathrm{ev}^F_{b,c}\big) \otimes \mathrm{colim}\big(\mathrm{ev}^F_{a,b}\big) & \xrightarrow[\mu]{} & \mathrm{colim}\big(\mathrm{ev}^F_{b,c} \boxtimes \mathrm{ev}^F_{a,b}\big) & \xrightarrow[\mathrm{colim}(\star,\theta)]{} & \mathrm{colim}\big(\mathrm{ev}^F_{a,c}\big)
\end{array}$$

*Unit map.* Let $\lambda$ denote the cocone of $\mathrm{colim}\big(\mathrm{ev}^F_{a,a}\big)$ and let $\eta$ be as in Proposition 3.8. The left square commutes, since the component of $F(i)$ in the cocone of $\mathrm{colim}(F)$ from Lemma 4.15 respects units. The right triangle commutes by definition of $\mathrm{colim}(\iota,\mathrm{id})$.

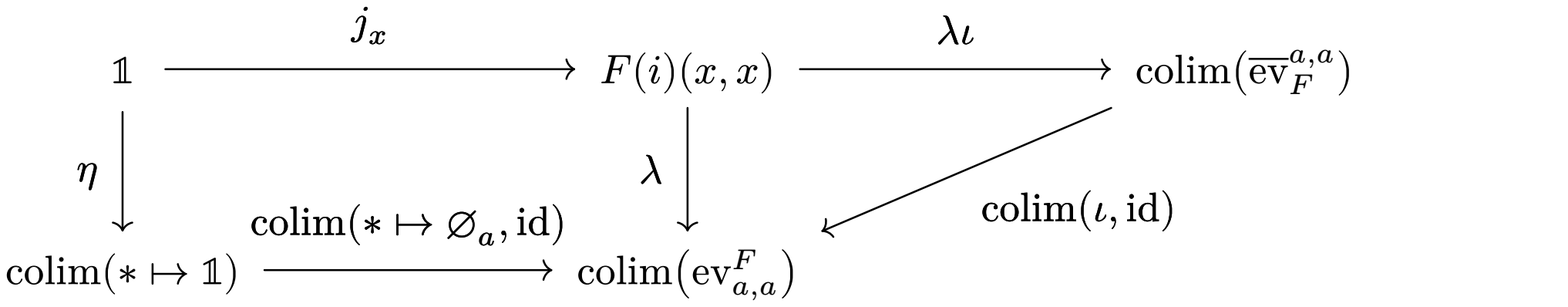

□

## 5. Lawvere metric spaces and ordinary categories

Definition 5.1 (Lawvere metric space) [Law73]. A **Lawvere metric space** is a category $M$ enriched in the monoidal poset $([0,\infty], \geq)$ with addition as the tensor product. Objects are *points*, hom-objects $M(a,b)$ are *distances* and composition identifies with the triangle inequality

$$M(b,c) + M(a,b) \geq M(a,c).$$

Remark 5.2. Lawvere metric spaces are strict generalizations of ordinary metric spaces, since they do not require symmetry $M(a,b) = M(b,a)$ and $M(a,b) = 0$ does not imply $a = b$.

A $[0,\infty]$-enriched functor $F\colon M \to N$ satisfies $FM(a,b) \geq N(Fa, Fb)$, i.e. it does not increase pairwise distances.

Lemma 5.3. Let $G\colon \mathcal{J} \to [0,\infty]$ be a diagram. Then $\inf_{j\in\mathcal{J}} G(j)$ is a colimit of $G$.

Note that $\operatorname{colim}(G)$ does not depend on morphisms in $\mathcal{J}$.

Corollary 5.4 (colimit of metric spaces). Let $F\colon \mathcal{I} \to \mathrm{Cat}_{[0,\infty]}$ be a diagram. Then the hom-distance $(\operatorname{colim} F)(a,b)$ is the infimum over sums of distances $\inf\left\{d_{F(i_n)}(x_n,y_n) + \dots + d_{F(i_1)}(x_1,y_1)\right\}$ with compatible endpoints $[y_{k-1}] = [x_k] \in \operatorname{colim}(\mathrm{Ob} \circ F)$ and $[x_1] = a$ and $[y_n] = b$.

Proof. It is an immediate consequence of Corollary 4.21 after applying Lemma 5.3 to $\operatorname{colim}_{\overline{\mathcal{S}}_F(a,b)}(\mathrm{ev}^F)$.□

Lemma 5.5 [Rie16, Theorem 3.5.11]. Assuming that the relevant coproducts and coequalizers exist, the colimit of $F\colon \mathcal{I} \to \mathcal{C}$ is isomorphic to the coequalizer of $d$ and $c$.

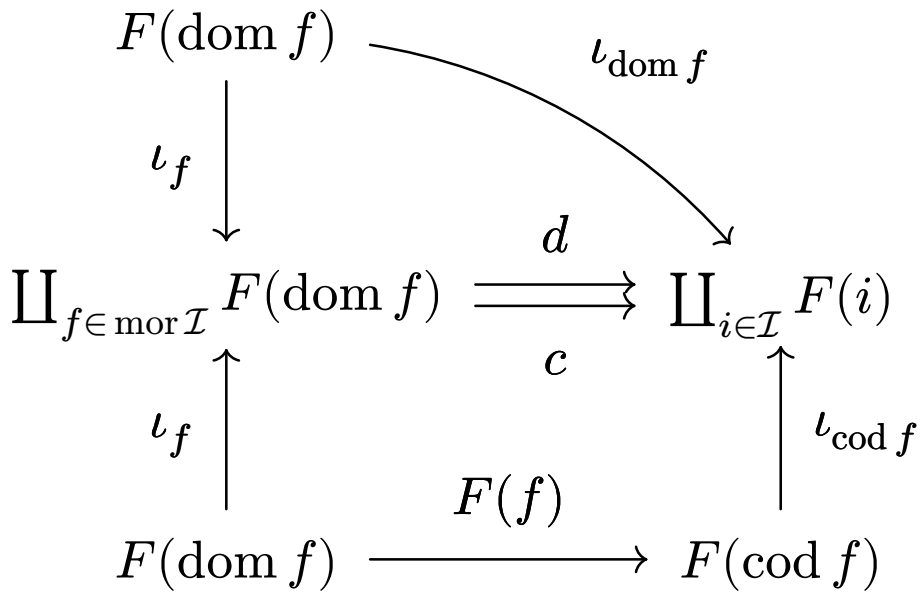

For pushouts of categories along embeddings, a closely related string description appears in [MS09, §2].

Corollary 5.6 (colimit of categories). Let $F\colon \mathcal{I} \to \mathrm{Cat}$ be a diagram. Then each morphism in $(\operatorname{colim} F)(a,b)$ between $a, b \in \operatorname{colim}(\mathrm{Ob} \circ F)$ is uniquely represented by an equivalence class of strings $[(i_n,\alpha_n), \dots, (i_1,\alpha_1)]$ where

$$\alpha_k \in \operatorname{mor} F(i_k), \quad \operatorname{cod}(\alpha_{k-1}) \sim \operatorname{dom}(\alpha_k), \quad \operatorname{cod}(\alpha_n) \in b, \quad \operatorname{dom}(\alpha_1) \in a$$

The equivalence relation of strings is generated by compositions of consecutive morphisms in the same $F(i)$ and edgewise transports along morphisms $f\colon i \to j$ in $\mathcal{I}$.

$$[\dots, (i,\beta), (i,\alpha), \dots] \sim [\dots, (i, \beta\circ\alpha), \dots]$$
$$[\dots, (i,\alpha), \dots] \sim [\dots, (j, (Ff)(\alpha)), \dots]$$

Proof. Interpret Cat as $\mathrm{Cat}_{\mathrm{Set}}$ with the cartesian monoidal structure $(\mathrm{Set}, \times, *)$. Apply Corollary 4.21 to $F$ and Lemma 5.5 to $\operatorname{colim}_{\overline{\mathcal{S}}_F(a,b)}(\mathrm{ev}^F)$. □

# 6. Rigidification of simplicial sets via necklaces

The category of simplicial sets sSet $:= [\Delta^{\mathsf{op}}, \mathrm{Set}]$ inherits the cartesian monoidal structure on Set: a tensor product $(X \times Y)_n := X_n \times Y_n$ and a unit $\Delta^0$, the terminal simplicial set.

Definition 6.1 (comma category) [Law63]. Given functors $F\colon \mathcal{I} \to \mathcal{C}$ and $G\colon \mathcal{J} \to \mathcal{C}$, the **comma category** $(F \downarrow G)$ has as objects triples $(i, j, Fi \to Gj)$ with $i \in \mathcal{I}$ and $j \in \mathcal{J}$ and as morphisms $(i, j, \alpha) \to (i', j', \alpha')$ pairs $(\beta\colon i \to i', \gamma\colon j \to j')$, such that the left diagram commutes in $\mathcal{C}$.

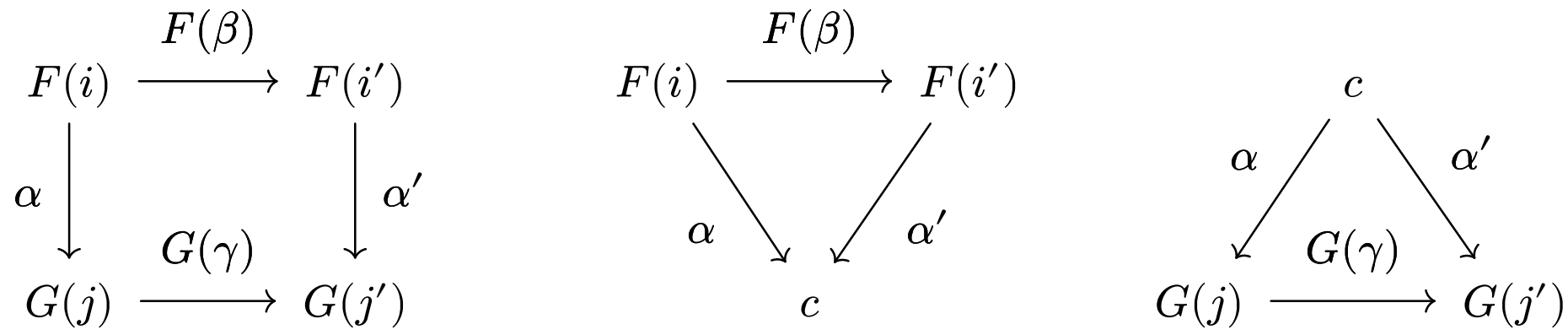

In the special case of $F\colon * \to \mathcal{C}$ with $* \mapsto c$ write $(F \downarrow G) =: (c \downarrow G)$ with objects $(j, c \to Gj)$. Analogously, $(F \downarrow c)$ has objects $(i, Fi \to c)$.

Remark 6.2. The comma category $(F \downarrow G)$ is a two-sided version of the category of elements (see Definition 4.1) of the profunctor $\mathcal{I}^{\mathsf{op}} \times \mathcal{J} \to \mathrm{Set}$ sending $(i, j) \mapsto \mathcal{C}(Fi, Gj)$.

Definition 6.3 (category of simplices). Let $\Delta$ denote the Yoneda embedding $\Delta\colon \Delta \to [\Delta^{\mathsf{op}}, \mathrm{Set}]$ given by $[n] \mapsto \Delta(-, [n]) =: \Delta^n$. Then the comma category $(\Delta \downarrow X)$ has as objects simplicial maps $\sigma\colon \Delta^n \to X$ and as morphisms commutative diagrams as on the right.

$$\begin{array}{ccc} \Delta^n & \searrow & \\ \downarrow & & X \\ \Delta^m & \nearrow & \end{array}$$

We defer the rigidification of standard simplices $\mathfrak{C}\Delta^n$ following [Lur08] to Section 6.1. We will also need to understand morphisms $\mathfrak{C}f$ given simplicial maps $f\colon \Delta^m \to \Delta^n$.

Definition 6.4 (rigidification) [Lur08, p. 29]. Let $X$ be a simplicial set. Define its **rigidification** $\mathfrak{C}X := \mathrm{colim}_{\Delta^n \to X}\, \mathfrak{C}\Delta^n$ as the colimit of the functor $F\colon (\Delta \downarrow X) \to \mathrm{Cat}_{\mathrm{sSet}}$ sending $\sigma\colon \Delta^n \to X$ to the simplex rigidification $\mathfrak{C}\Delta^n$ defined below in Definition 6.7.

## 6.1. Rigidification of standard simplex

Lemma 6.5. The nerve $N\colon \mathrm{Cat} \to \mathrm{sSet}$ given by $(N\mathcal{C})_n := [[n], \mathcal{C}]$ is a strong monoidal functor between cartesian monoidal categories $(\mathrm{Cat}, \times, *)$ and $(\mathrm{sSet}, \times, \Delta^0)$.

Proof. $\Delta^0 \cong N(*)$ is obvious. Continuity of $N$ yields a natural isomorphism $N\mathcal{C} \times N\mathcal{D} \cong N(\mathcal{C} \times \mathcal{D})$. □

Definition 6.6. Define a Cat-enriched category $P$ with objects $\mathbb{N}_0$. For $i \leq j$ let $P(i, j)$ be the poset category of subsets of $\{i, ..., j\}$ containing $\{i, j\}$ and set $P(i, j) = \varnothing$ for $i > j$. Composition is given by union

$$\cup : P(j, k) \times P(i, j) \to P(i, k).$$

This gives rise to a cosimplicial object in $\mathrm{Cat}_{\mathrm{Cat}}$, i.e. a functor $P\colon \Delta \to \mathrm{Cat}_{\mathrm{Cat}}$, by taking $P[n]$ to be the restriction of $P$ to objects $\{0, ..., n\}$. Given a monotone map $f\colon [m] \to [n]$ in $\Delta$, define the Cat-enriched functor $P(f)\colon P[m] \to P[n]$ on objects by $i \mapsto f(i)$ and on hom-categories $P(i, j) \to P(fi, fj)$ by $I \mapsto f(I)$.

By Lemma 2.5, the nerve $N$ induces a functor $N_*\colon \mathrm{Cat}_{\mathrm{Cat}} \to \mathrm{Cat}_{\mathrm{sSet}}$ preserving tensor and unit up to coherent isomorphism.

Definition 6.7 (simplex rigidification) [Lur08, Definition 1.1.5.1, Definition 1.1.5.3]. Let $\mathfrak{C}$ denote the composition $N_* \circ P\colon \Delta \to \mathrm{Cat}_{\mathrm{Cat}} \to \mathrm{Cat}_{\mathrm{sSet}}$ with $P$ as in Definition 6.6. We call $\mathfrak{C}\Delta^n := N_*P[n]$ the **rigidification** of $\Delta^n$.

In Appendix C, we discuss an alternative description of $\mathfrak{C}\Delta^n$, which might provide deeper intuition. However, Definition 6.7 is more suitable for the computations in this section. In Lemma C.5, we prove that both simplex rigidifications $\Delta \to \mathsf{Cat}_{\mathrm{sSet}}$ are naturally isomorphic.

## 6.2. FROM STRINGS TO NECKLACES

DESCRIPTION 6.8. We unravel the category of non-empty strings $\overline{\mathcal{S}}_F$ from Definition 4.17 for $F$ as in Definition 6.4. Edges $(\sigma, x, y)$ are given by $\sigma\colon \Delta^n \to X$ and $x, y \in \mathrm{Ob}(\mathfrak{C}\Delta^n) = \mathrm{Ob}(\Delta^n)$. A string $[(\sigma_m, x_m, y_m), ..., (\sigma_1, x_1, y_1)]$ has compatible endpoints $[y_{k-1}] = [x_k] \in \mathrm{colim}(\mathrm{Ob} \circ F) = X_0$, hence $\sigma_{k-1}(y_{k-1}) = \sigma_k(x_k) \in X_0$. Morphisms are generated by transports $\Delta^n \to \Delta^m$ over $X$ and contractions $[..., (\sigma, y, z), (\sigma, x, y), ...] \to [..., (\sigma, x, z), ...]$.

DEFINITION 6.9 (necklace) [DS11, p. 6]. A **necklace** is a simplicial set $T := \Delta^{n_1} \vee \cdots \vee \Delta^{n_m}$, where the final vertex of the **bead** $\Delta^{n_k}$ is glued to the initial vertex of $\Delta^{n_{k+1}}$. Let $J_T$ denote the **joints** and $V_T$ the **vertices**. We allow degenerate beads $\Delta^0$. The **category of necklaces** Nec has as morphisms endpoint-preserving simplicial maps. Define the comma category $(\mathrm{Nec} \downarrow X) := ((\mathrm{Nec} \to \mathrm{sSet}) \downarrow X)$ and the restriction $(\mathrm{Nec} \downarrow X)_{a,b}$ to maps $T \to X$ sending the initial and final necklace vertices to $a$ and $b$ respectively.

We need to understand endpoint-preserving maps between necklaces, since they generate an equivalence relation in Corollary 6.18. We have to prove that the equivalence classes agree with the hom-objects in Theorem 4.14.

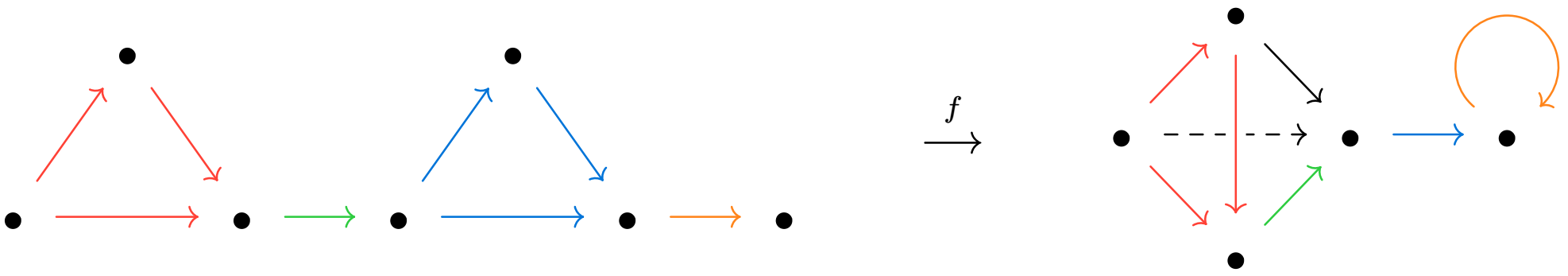

LEMMA 6.10. Let $f\colon T \to U$ be an endpoint-preserving map of necklaces. Then for each bead $\Delta^n \subseteq T$, there exists a bead of $U$ containing $f(\Delta^n)$ and $J_U \subseteq f(J_T)$.

PROOF. Assume $y \in J_U - f(J_T) \neq \varnothing$. There exists a bead $\Delta^n$ of $T$ with endpoints $x, z \in J_T$, such that $y$ lies between $f(x)$ and $f(z)$. However, $U$ cannot contain the 1-simplex $\Delta^1 \subseteq f(\Delta^n)$ from $f(x)$ to $f(z)$. □

DEFINITION 6.11 (category of directed, non-empty strings). For $F\colon (\Delta \downarrow X) \to \mathsf{Cat}_{\mathrm{sSet}}$ as in Definition 6.4, define the full subcategory $\vec{\mathcal{S}}_X(a,b) \subset \overline{\mathcal{S}}_F(a,b)$ of **directed strings** that consist of edges $(\sigma, x, y)$ satisfying $x \leq y$.

We restrict to directed strings for well-definedness of the following "necklace replacements".

DEFINITION 6.12 (necklace replacement). Let $[(\sigma_m, x_m, y_m), ..., (\sigma_1, x_1, y_1)] \in \vec{\mathcal{S}}_X(a,b)$. Define its **necklace replacement** as the necklace $T := \Delta^{y_1 - x_1} \vee \cdots \vee \Delta^{y_m - x_m}$ together with the map $T \to X$, defined bead-wise as

$$\Delta^{y_k - x_k} \xhookrightarrow{+x_k} \Delta^{n_k} \xrightarrow{\sigma_k} X.$$

The joints $J_T$ of $T$ map to $a = [x_1], ..., [y_{k-1}] = [x_k], ..., [y_m] = b \in X_0$. Taking necklace replacements extends to a functor $\mathrm{repl}_X^{a,b}\colon \vec{\mathcal{S}}_X(a,b) \to (\mathrm{Nec} \downarrow X)_{a,b}$. Edgewise transports $f\colon \Delta^n \to \Delta^m$ over $X$ restrict to $\Delta^{y-x} \to \Delta^{fy-fx}$ and contractions are sent to $\Delta^{y-x} \vee \Delta^{z-y} \hookrightarrow \Delta^{z-x}$.

LEMMA 6.13. The necklace replacement functor $\mathrm{repl}_X$ from Definition 6.12 is final, i.e. restricting diagrams along $\mathrm{repl}_X$ preserves their colimit.

PROOF. Abbreviate $\mathrm{repl} := \mathrm{repl}_X$ and suppress superscript endpoints. By [Mac98, §IX.3], we need to prove that for each simplicial map $t\colon T \to X$, the comma category $(t \downarrow \mathrm{repl})$ is non-empty and connected.

- Objects are pairs $(s, f)$ of non-empty strings $s \in \vec{\mathcal{S}}_X$ and endpoint-preserving necklace maps $f\colon t \to \mathrm{repl}(s)$ over $X$.
- A morphism $(s, f) \to (s', f')$ is a string morphism $\varphi\colon s \to s'$, such that $f' = \mathrm{repl}(\varphi) \circ f$.

*Non-empty.* Write $T \coloneqq \Delta^{n_1} \vee \cdots \vee \Delta^{n_m}$ and note that each $t|_{\Delta^{n_i}} \in (\Delta \downarrow X)$. The string $s \coloneqq [(t|_{\Delta^{n_m}}, 0, n_m), ..., (t|_{\Delta^{n_1}}, 0, n_1)]$ maps to $\mathrm{repl}(s) = t$.

*Connected.* It suffices to show that any necklace morphism $f\colon t \to \mathrm{repl}(s')$ over $X$ with $\mathrm{repl}(s')\colon T' \to X$ is connected to $\mathrm{id}\colon t \to \mathrm{repl}(s)$ by a zigzag. We can assume that $s'$ either consists of a single degenerate edge, or contains no degenerate edges $(\sigma, x, x)$, by contracting after a zigzag of transports as in the left diagram. This lifts to a zigzag in $(t \downarrow \mathrm{repl})$.

$$[..., (\tau, y, z), (\sigma, x, x), ...] \leftrightsquigarrow [..., (\tau, y, z), (\tau, y, y), ...] \to [..., (\tau, y, z), ...]$$

We will construct a string morphism $\varphi\colon s \to s'$ satisfying $f = \mathrm{repl}(\varphi)\colon T \to T'$ as in the right diagram.

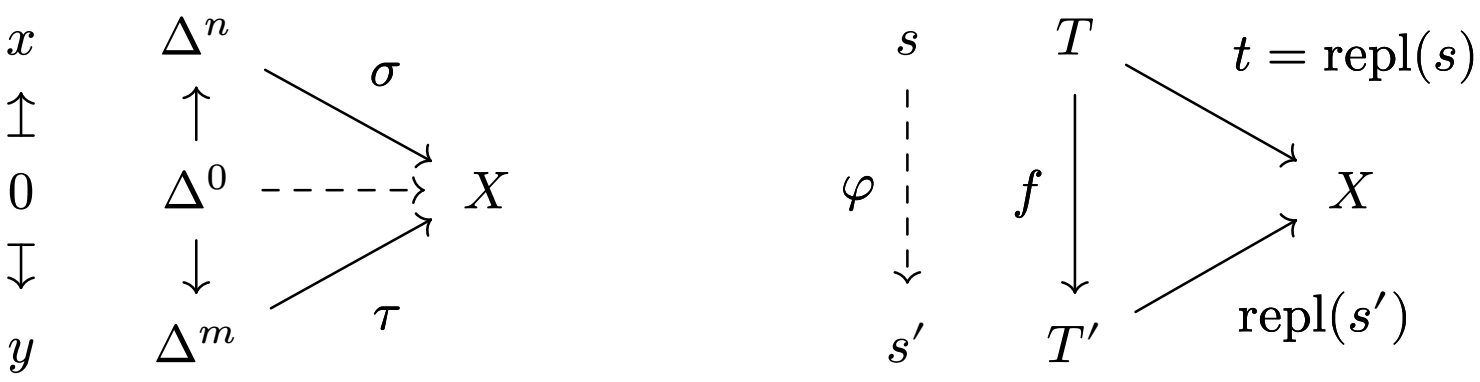

Write $s' \coloneqq [(\sigma_\ell, x_\ell, y_\ell), ..., (\sigma_1, x_1, y_1)]$ with $x_i < y_i$. By Lemma 6.10, $J_{T'} \subseteq f(J_T)$. Thus, each restriction $f|_{\Delta^{n_i}}$ lands in a bead of $T'$ and underlies a morphism in $(\Delta \downarrow X)$. Define $\varphi$ to be the corresponding edgewise transports $(t|_{\Delta^{n_i}}, 0, n_i) \to (\sigma_j, u_i, v_i)$ with $x_j \le u_i \le v_i \le y_j$, followed by contractions of all consecutive edges in the same $\sigma_j$. □

## 6.3. From string evaluation to flags of necklaces

Lemma 6.14. For $F\colon (\Delta \downarrow X) \to \mathrm{Cat}_{\mathrm{sSet}}$ as in Definition 6.4, denote by $\vec{\mathrm{ev}}_X^{a,b}$ the restriction of string evaluation (see Definition 4.8 and Corollary 4.21) $\overline{\mathrm{ev}}_F^{a,b}$ to the full subcategory of directed strings $\vec{\mathcal{S}}_X(a,b)$ from Definition 6.11. The subcategory inclusion induces an isomorphism $\mathrm{colim}\big(\vec{\mathrm{ev}}_X^{a,b}\big) \cong \mathrm{colim}\big(\overline{\mathrm{ev}}_F^{a,b}\big)$.

Proof. Note that there exist no morphisms from directed to undirected $F$-strings. The functor $F$ sends a simplicial map $\sigma\colon \Delta^n \to X$ to the rigidification $\mathfrak{C}\Delta^n = N_* P[n]$. By Definition 6.6, $P(x,y) = \varnothing$ for $x > y$, so a string $s \in \overline{\mathcal{S}}_F(a,b)$ containing the edge $(\sigma, x, y)$ evaluates to

$$\overline{\mathrm{ev}}(s) = \cdots \times NP(x,y) \times \cdots = \varnothing.$$

Via the unique universal morphism from this initial object in sSet, any cocone $\vec{\mathrm{ev}}_X^{a,b} \Rightarrow \Delta K$ extends to a cocone $\overline{\mathrm{ev}}_F^{a,b} \Rightarrow \Delta K$. As in the proof of Lemma 4.20, this extension is inverse to the restriction induced by $\vec{\mathcal{S}}_X(a,b) \subset \overline{\mathcal{S}}_F(a,b)$. By naturality in $K$, Lemma 4.18 concludes the proof. □

Definition 6.15. Define the functor flag: $(\mathrm{Nec} \downarrow X) \to \mathrm{sSet}$ on objects $T \to X$ by $n$-simplices

$$\mathrm{flag}(T \to X)_n \coloneqq \{J_T \subseteq T^0 \subseteq \cdots \subseteq T^n \subseteq V_T\}$$

and on endpoint-preserving maps $f\colon T \to U$ over $X$ by $(T^0 \subseteq \cdots \subseteq T^n) \mapsto (f(T^0) \subseteq \cdots \subseteq f(T^n))$. The degeneracy map $s_i$ duplicates $T^i$, while the face map $d_i$ omits $T^i$ from the representative.

By Lemma 6.10, the functor is well-defined on morphisms.

Lemma 6.16. Let $F$ be as in Definition 6.4. There exists a natural isomorphism $\vec{\mathrm{ev}}_X^{a,b} \simeq \mathrm{flag} \circ \mathrm{repl}_X^{a,b}$.

Proof. We start with $\vec{\mathrm{ev}} \coloneqq \vec{\mathrm{ev}}_X^{a,b}$ and compose multiple natural isomorphisms $\eta_3 \circ \eta_2 \circ \eta_1$ to obtain $\mathrm{flag} \circ \mathrm{repl}$.

*Strings.* Consider an $F$-string $s \coloneqq [(\sigma_m, x_m, y_m), ..., (\sigma_1, x_1, y_1)]$. By symmetry of the cartesian monoidal structure on sSet, we can reverse the order of the product up to a natural isomorphism $\eta_1$. Use continuity of the nerve $N$ to obtain the second natural isomorphism $\eta_2$.

$$\begin{aligned}\vec{\mathrm{ev}}(s) &= \mathfrak{C}\Delta^{n_m}(x_m, y_m) \times \cdots \times \mathfrak{C}\Delta^{n_1}(x_1, y_1) \\ &\overset{\eta_1}{\cong} NP(x_1, y_1) \times \cdots \times NP(x_m, y_m) \\ &\overset{\eta_2}{\cong} N(P(x_1, y_1) \times \cdots \times P(x_m, y_m))\end{aligned} \tag{3}$$

Let $T := \mathrm{repl}(s)$. Since endpoints coincide, i.e. $\sigma_{k-1}(y_{k-1}) = \sigma_k(x_k)$, a tuple $(I_1, ..., I_m) \in \bigtimes_k P(x_k, y_k)$ is naturally identified with a subset of $V_T$ containing $J_T$. Thus, the entire nerve $N \bigtimes_k P(x_k, y_k)$ is identified with $\mathrm{flag}(T)$ via a natural isomorphism $\eta_3$.

*Edgewise transport.* A morphism $f\colon (\Delta^m \to X, x, y) \to (\Delta^n \to X, x', y')$ between edges comes from an underlying morphism $[m] \to [n]$ in $\Delta$ with $f(x) = x'$ and $f(y) = y'$. Recall from Definition 6.7 that $F(f) = N_* P(f)$ and the simplicial map $N(Pf)_{x,y}\colon NP(x, y) \to NP(fx, fy)$ is given on $k$-simplices by

$$(I_0 \subseteq \cdots \subseteq I_k) \mapsto (f(I_0) \subseteq \cdots \subseteq f(I_k)).$$

An edgewise transport $[f_m, ..., f_1]$ is sent to $Ff_m \times \cdots \times Ff_1$ by $\vec{\mathrm{ev}}$. By naturality of $\eta_2$ in Equation 3, the nerve commutes with products of functors. To reduce visual noise, we only track one edge $(\sigma \circ \Delta f, x, y)$, where $\Delta$ is the Yoneda embedding from Definition 6.3. Keep in mind that $\eta_1$ reverses the order of the product on the right.

$$\begin{array}{ccc} \vec{\mathrm{ev}}[..., (\sigma \circ \Delta f, x, y), ...] & \xrightarrow{\cong} & N(\cdots \times P(x, y) \times \cdots) \\ {\scriptstyle \vec{\mathrm{ev}}[..., f, ...]}\downarrow & & \downarrow{\scriptstyle N(\cdots \times Pf \times \cdots)} \\ \vec{\mathrm{ev}}[..., (\sigma, fx, fy), ...] & \xrightarrow{\cong} & N(\cdots \times P(fx, fy) \times \cdots) \end{array}$$

Since each edge morphism $f_k\colon (\sigma_k, x_k, y_k) \to (\sigma'_k, x'_k, y'_k)$ satisfies $f_k(x_k) = x'_k$ and $f_k(y_k) = y'_k$, the consecutive endpoint identifications by $\eta_3$ factor through edgewise transports and $N(Pf_1 \times \cdots \times Pf_m)$ descends to $\mathrm{flag}(\mathrm{repl}[f_m, ..., f_1])$.

$$\sigma'_{k-1}(y'_{k-1}) = \sigma_{k-1}(y_{k-1}) = \sigma_k(x_k) = \sigma'_k(x'_k)$$

*Contractions.* The string evaluation functor $\vec{\mathrm{ev}}$ sends contractions to

$$\cdots \times NP(y, z) \times NP(x, y) \times \cdots \cong \cdots \times N(P(y, z) \times P(x, y)) \times \cdots \xrightarrow{\cdots \times N(\cup) \times \cdots} \cdots \times NP(x, z) \times \cdots$$

The natural isomorphisms $\eta_1$ and $\eta_2$ in Equation 3 identify this contraction with

$$N(\cdots \times P(x, y) \times P(y, z) \times \cdots) \xrightarrow{N(\cdots \times \cup \times \cdots)} N(\cdots \times P(x, z) \times \cdots)$$

Remember that the necklace replacement functor maps the contraction to the inclusion $\iota\colon \Delta^{y-x} \vee \Delta^{z-y} \hookrightarrow \Delta^{z-x}$. After $\eta_3$ identifies endpoints in $P(x, y) \times P(y, z)$, the map $N(\cdots \times \cup \times \cdots)$ becomes $\mathrm{flag}(\cdots \vee \iota \vee \cdots)$. □

Theorem 6.17 [DS11, Proposition 3.8, Proposition 4.3]. The hom-sSet $\mathfrak{C}X(a, b)$ is isomorphic to the colimit of the restriction to $(\mathrm{Nec} \downarrow X)_{a,b}$ of flag as in Definition 6.15.

Proof. Set $F\colon (\Delta \downarrow X) \to \mathrm{Cat}_{\mathrm{sSet}}$ as in Definition 6.4. By Corollary 4.21, the rigidification $\mathfrak{C}X = \mathrm{colim}\, F$ has as hom-sSet $\mathfrak{C}X(a, b) \cong \mathrm{colim}\left(\overline{\mathrm{ev}}_F^{a,b}\right)$. By Lemma 6.14, restricting to directed strings preserves the colimit. By Lemma 6.16, $\vec{\mathrm{ev}}_X^{a,b} \simeq \mathrm{flag} \circ \mathrm{repl}_X^{a,b}$ are naturally isomorphic, so their colimits coincide by Proposition 3.3. By Lemma 6.13, $\mathrm{repl}_X$ is final, so restricting flag along $\mathrm{repl}_X$ preserves the colimit.

$$\mathrm{colim}\, \vec{\mathrm{ev}}_X^{a,b} \cong \operatorname*{colim}_{\vec{\mathcal{S}}(a,b)}(\mathrm{flag} \circ \mathrm{repl}_X) \cong \operatorname*{colim}_{(\mathrm{Nec} \downarrow X)_{a,b}} \mathrm{flag}$$

□

COROLLARY 6.18 [DS11, Corollary 4.4]. Let $X \in \mathrm{sSet}$ and $a, b \in X_0$. An $n$-simplex in $\mathfrak{C}X(a,b)$ is uniquely represented by an equivalence class $\left[T \to X, \vec{T}\right]$, where

- $T \to X$ is a simplicial map sending $\alpha_T \mapsto a$ and $\omega_T \mapsto b$.
- $\vec{T}$ is a flag $J_T \subseteq T^0 \subseteq \cdots \subseteq T^n \subseteq V_T$.

The equivalence relation is generated by $\left(T \to X, \vec{T}\right) \sim \left(U \to X, \vec{U}\right)$ if there exists $f\colon T \to U$ over $X$ with $\vec{U} = f_* \vec{T}$. The degeneracy map $s_i$ duplicates $T^i$, while the face map $d_i$ omits $T^i$ from the representative.

PROOF. By Lemma 5.5, the colimit of flag: $(\mathrm{Nec} \downarrow X)_{a,b} \to \mathrm{sSet}$ is isomorphic to the coequalizer of $d$ and $c$, where $T \to U \to X$ denotes an endpoint-preserving necklace morphism $T \to U$ over $X$.

$$\coprod_{T \to U \to X} \mathrm{flag}(T \to X) \underset{c}{\overset{d}{\rightrightarrows}} \coprod_{T' \to X} \mathrm{flag}(T' \to X)$$

□

# A. COHERENCE VERIFICATIONS

LEMMA A.1. The monoidal structure of $\int[-, \mathcal{V}]$ from Lemma 3.5 satisfies the triangle and pentagon identities.

PROOF. *Triangle identity.* Write $\mathbb{1}_* := (* \mapsto \mathbb{1})$.

$$\begin{array}{ccc} ((\mathcal{I} \times *) \times \mathcal{J}, (F \boxtimes \mathbb{1}_*) \boxtimes G) & \xrightarrow{\left(\alpha^{\mathrm{Cat}}_{\mathcal{I},*,\mathcal{J}}, \alpha^{\mathcal{V}}_{F,\mathbb{1}_*,G}\right)} & (\mathcal{I} \times (* \times \mathcal{J}), F \boxtimes (\mathbb{1}_* \boxtimes G)) \\ & {\scriptstyle \left(r^{\mathrm{Cat}}_{\mathcal{I}}, r^{\mathcal{V}}_F\right) \boxtimes (\mathrm{id}_{\mathcal{J}}, \mathrm{id}_G)} \searrow \quad \swarrow {\scriptstyle (\mathrm{id}_{\mathcal{I}}, \mathrm{id}_F) \boxtimes \left(l^{\mathrm{Cat}}_{\mathcal{J}}, l^{\mathcal{V}}_G\right)} & \\ & (\mathcal{I} \times \mathcal{J}, F \boxtimes G) & \end{array}$$

Recalling composition in $\int[-, \mathcal{V}]$, we spell out the diagram in $[(\mathcal{I} \times *) \times \mathcal{J}, \mathcal{V}]$. By the triangle identity in Cat, we have equality at the bottom.

$$\begin{array}{ccc} (F \boxtimes \mathbb{1}_*) \boxtimes G & \overset{\alpha^{\mathcal{V}}_{F,\mathbb{1}_*,G}}{\Longrightarrow} & (F \boxtimes (\mathbb{1}_* \boxtimes G)) \circ \alpha^{\mathrm{Cat}}_{\mathcal{I},*,\mathcal{J}} \\ {\scriptstyle r^{\mathcal{V}}_F \boxtimes \mathrm{id}_G} \Downarrow & & \Downarrow {\scriptstyle \left(\mathrm{id}_F \boxtimes l^{\mathcal{V}}_G\right)\alpha^{\mathrm{Cat}}_{\mathcal{I},*,\mathcal{J}}} \\ (F \boxtimes G) \circ \left(r^{\mathrm{Cat}}_{\mathcal{I}} \times \mathrm{id}_{\mathcal{J}}\right) & = & (F \boxtimes G) \circ \left(\mathrm{id}_{\mathcal{I}} \times l^{\mathrm{Cat}}_{\mathcal{J}}\right) \circ \alpha^{\mathrm{Cat}}_{\mathcal{I},*,\mathcal{J}} \end{array}$$

Insert some $((i, *), j)$ and use the triangle identity in $\mathcal{V}$.

$$\begin{array}{ccc} (F(i) \otimes \mathbb{1}) \otimes G(j) & \xrightarrow{\alpha^{\mathcal{V}}_{F(i),\mathbb{1},G(j)}} & F(i) \otimes (\mathbb{1} \otimes G(j)) \\ & {\scriptstyle r^{\mathcal{V}}_{F(i)} \boxtimes \mathrm{id}} \searrow \quad \swarrow {\scriptstyle \mathrm{id} \boxtimes l^{\mathcal{V}}_{G(j)}} & \\ & F(i) \otimes G(j) & \end{array}$$

*Pentagon identity.* Analogously, use the pentagon identity of Cat and $\mathcal{V}$. □

LEMMA A.2. The colimit functor $\mathsf{colim}\colon \int[-, \mathcal{V}] \to \mathcal{V}$ from Proposition 3.8 satisfies unital and associative coherence.

PROOF. Let $\iota$ denote the colimiting cocone of $\mathsf{colim}(F \boxtimes F')$. We have already seen that $\iota = \Delta\mu \circ \kappa \boxtimes \kappa'$.

*Unitality.* We only prove left unitality, since right unitality is shown analogously. Let $\kappa$ denote the cocone of $\mathsf{colim}(F)$. Note that the components $\mathrm{id} \otimes \kappa_i\colon \mathbb{1} \otimes F(i) \to \mathbb{1} \otimes \mathsf{colim}\, F$ define a colimiting cocone of $\mathbb{1} \otimes \mathsf{colim}(F) \cong \mathsf{colim}(\mathbb{1} \otimes F)$. By the colimit's universal property, it suffices to show that both morphisms $\mathbb{1} \otimes \mathsf{colim}_{\mathcal{I}} F \to \mathsf{colim}_{\mathcal{I}} F$ in the following diagram coincide, after precomposing with any $\mathrm{id}_{\mathbb{1}} \otimes \kappa_i$.

$$\begin{array}{ccc} \mathbb{1}\otimes \operatorname{colim}(\mathcal{I},F) & \xrightarrow{\eta\otimes \mathrm{id}} & \operatorname{colim}(*\mapsto \mathbb{1})\otimes\operatorname{colim}(\mathcal{I},F) \\ {\scriptstyle l}\downarrow & & \downarrow{\scriptstyle \mu} \\ \operatorname{colim}(\mathcal{I},F) & \xleftarrow{\operatorname{colim}(l^{\mathrm{Cat}},l)} & \operatorname{colim}((*\mapsto\mathbb{1})\boxtimes(\mathcal{I},F)) \end{array}$$

Let $\iota$ denote the colimiting cocone of $\operatorname{colim}((*\mapsto\mathbb{1})\boxtimes F)$. In the following diagram, the upper triangle commutes by functoriality of $\otimes$, the left square commutes by naturality of $l$, the right triangle commutes by definition of $\mu$ and the middle square commutes by definition of $\operatorname{colim}(l^{\mathrm{Cat}},l)$.

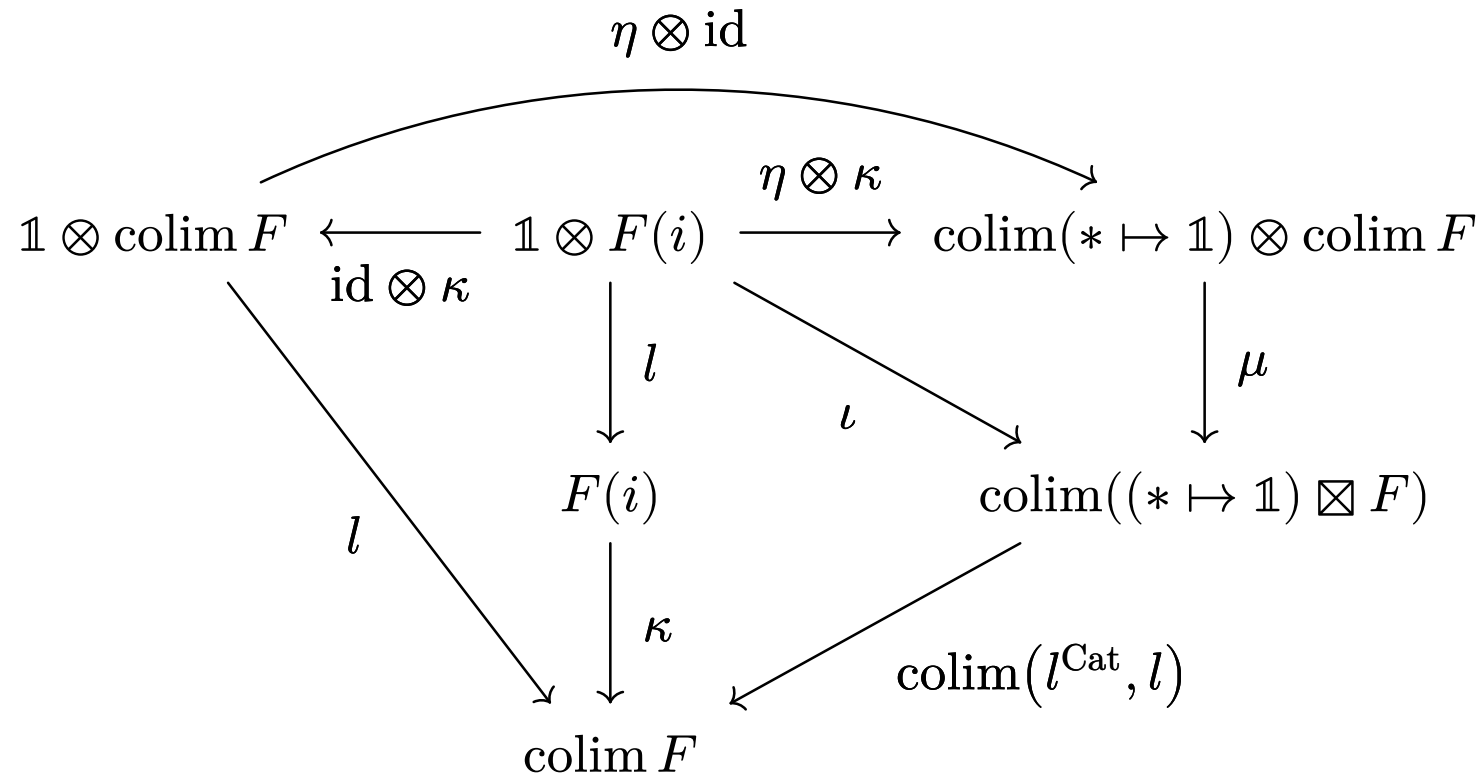

*Associativity.* Let $\kappa$, $\kappa'$, $\kappa''$, $\iota$ and $\iota'$ denote the colimiting cocones of the colimits of $F$, $G$, $H$, $F\boxtimes G$ and $G\boxtimes H$ respectively. We prove associativity coherence for the inverse isomorphism $\mu^{-1}$.

$$\begin{array}{ccc} \operatorname{colim}_{(\mathcal{I}\times\mathcal{J})\times\mathcal{K}}((F\boxtimes G)\boxtimes H) & \xrightarrow{\operatorname{colim}(\alpha^{\mathrm{Cat}},\alpha)} & \operatorname{colim}_{\mathcal{I}\times(\mathcal{J}\times\mathcal{K})}(F\boxtimes(G\boxtimes H)) \\ {\scriptstyle \mu^{-1}}\downarrow & & \downarrow{\scriptstyle \mu^{-1}} \\ \operatorname{colim}_{\mathcal{I}\times\mathcal{J}}(F\boxtimes G)\otimes\operatorname{colim}_{\mathcal{K}} H & & \operatorname{colim}_{\mathcal{I}} F\otimes\operatorname{colim}_{\mathcal{J}\times\mathcal{K}}(G\boxtimes H) \\ {\scriptstyle \mu^{-1}\otimes\mathrm{id}}\downarrow & & \downarrow{\scriptstyle \mathrm{id}\otimes\mu^{-1}} \\ (\operatorname{colim}_{\mathcal{I}} F\otimes\operatorname{colim}_{\mathcal{J}} G)\otimes\operatorname{colim}_{\mathcal{K}} H & \xrightarrow{\alpha} & \operatorname{colim}_{\mathcal{I}} F\otimes(\operatorname{colim}_{\mathcal{J}} G\otimes\operatorname{colim}_{\mathcal{K}} H) \end{array}$$

By the universal property of $\operatorname{colim}((F\boxtimes G)\boxtimes H)$, it suffices to show commutativity in $[(\mathcal{I}\times\mathcal{J})\times\mathcal{K},\mathcal{V}]$ after precomposition with its colimiting cocone. With $\Delta=\Delta_{(\mathcal{I}\times\mathcal{J})\times\mathcal{K}}$, the outer triangles commute by definition of $\mu^{-1}$ and the middle square commutes by naturality of $\alpha$.

$$\begin{array}{ccccc} & (F\boxtimes G)\boxtimes H & \xRightarrow{\alpha_{F,G,H}} & F\boxtimes(G\boxtimes H)\circ\alpha^{\mathrm{Cat}}_{\mathcal{I},\mathcal{J},\mathcal{K}} & \\ {\scriptstyle \iota\boxtimes\kappa''}\swarrow & \Downarrow{\scriptstyle (\kappa\boxtimes\kappa')\boxtimes\kappa''} & & \Downarrow{\scriptstyle (\kappa\boxtimes(\kappa'\boxtimes\kappa''))\alpha^{\mathrm{Cat}}} & \searrow{\scriptstyle (\kappa\boxtimes\iota')\alpha^{\mathrm{Cat}}} \\ \Delta(\operatorname{colim}(F\boxtimes G)\otimes\operatorname{colim} H) & & & & \Delta(\operatorname{colim} F\otimes\operatorname{colim}(G\boxtimes H)) \\ {\scriptstyle \Delta(\mu^{-1}\otimes\mathrm{id})}\searrow & & & & \swarrow{\scriptstyle \Delta(\mathrm{id}\otimes\mu^{-1})} \\ & \Delta((\operatorname{colim} F\otimes\operatorname{colim} G)\otimes\operatorname{colim} H) & \xRightarrow{\Delta\alpha} & \Delta(\operatorname{colim} F\otimes(\operatorname{colim} G\otimes\operatorname{colim} H)) & \end{array}$$

□

Lemma A.3. $\mathbb{S}_F$ from Definition 4.11 satisfies the axioms of a $\int[-,\mathcal{V}]$-enriched category.

Proof. Untangle the compositions in $\int[-,\mathcal{V}]$ as in the proof of Lemma A.1.

*Associativity.* We will prove commutativity of this diagram in $\int[-,\mathcal{V}]$.

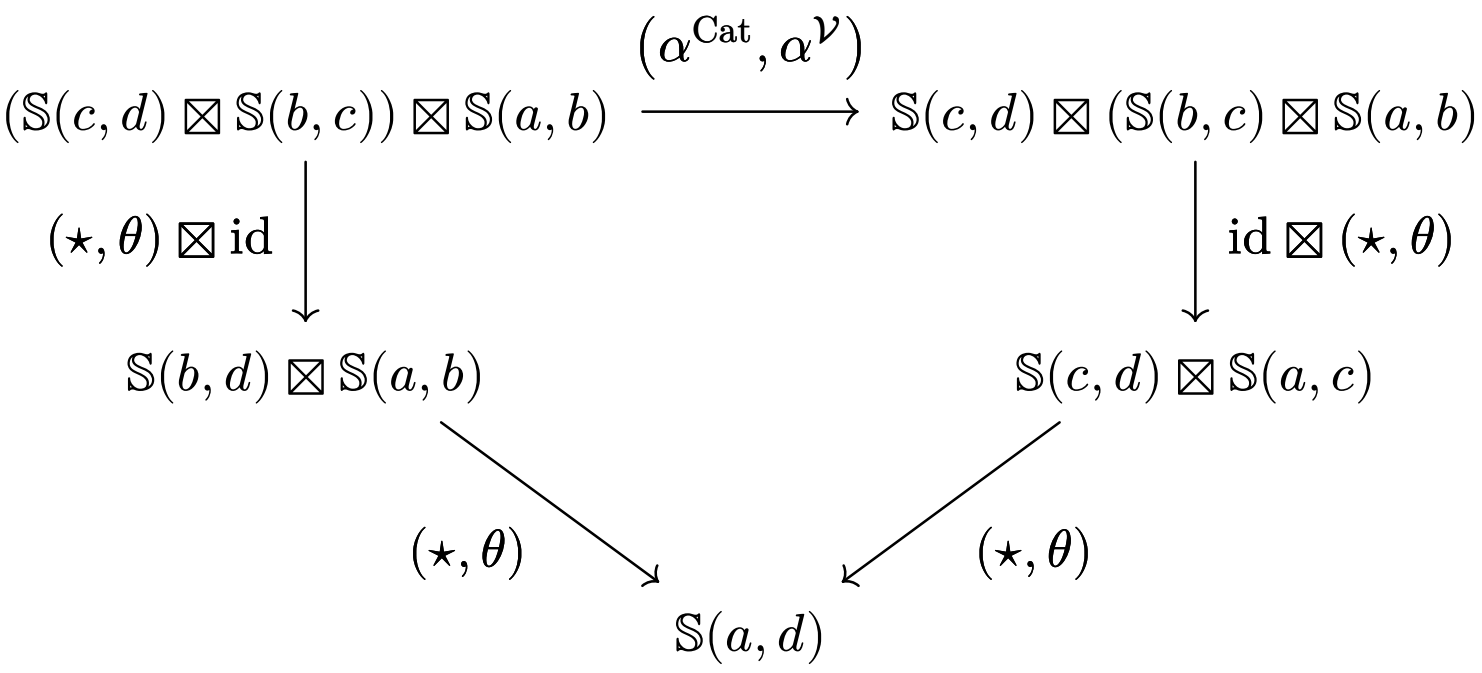

We have equality at the bottom of the following diagram in $[(\mathcal{S}(c,d)\times\mathcal{S}(b,c))\times\mathcal{S}(a,b),\mathcal{V}]$ by associativity for the Cat-enriched $\mathcal{S}_F$.

$$\begin{array}{ccc}
\left(\mathrm{ev}^{c,d}\boxtimes\mathrm{ev}^{b,c}\right)\boxtimes\mathrm{ev}^{a,b} & \overset{\alpha^{\mathcal{V}}_{\mathrm{ev,ev,ev}}}{\Longrightarrow} & \mathrm{ev}^{c,d}\boxtimes\left(\mathrm{ev}^{b,c}\boxtimes\mathrm{ev}^{a,b}\right)\circ\alpha^{\mathrm{Cat}} \\
\theta\boxtimes\mathrm{id}\Downarrow & & \Downarrow(\mathrm{id}\boxtimes\theta)\alpha^{\mathrm{Cat}} \\
\mathrm{ev}^{b,d}\left(-\star_{b,c,d}-\right)\boxtimes\mathrm{ev}^{a,b} & & \left(\mathrm{ev}^{c,d}\boxtimes\mathrm{ev}^{a,c}\left(-\star_{a,b,c}-\right)\right)\circ\alpha^{\mathrm{Cat}} \\
\theta\Downarrow & & \Downarrow\theta\alpha^{\mathrm{Cat}} \\
\mathrm{ev}^{a,d}\left(\left(-\star_{b,c,d}-\right)\star_{a,b,d}-\right) & = & \mathrm{ev}^{a,d}\left(-\star_{a,c,d}\left(-\star_{a,b,c}-\right)\right)\circ\alpha^{\mathrm{Cat}}
\end{array}$$

Insert some $((s_3,s_2),s_1)$. The pentagon commutes by definition of $\theta$.

$$\begin{array}{ccc}
(\mathrm{ev}(s_3)\otimes\mathrm{ev}(s_2))\otimes\mathrm{ev}(s_1) & \xrightarrow{\alpha^{\mathcal{V}}} & \mathrm{ev}(s_3)\otimes(\mathrm{ev}(s_2)\otimes\mathrm{ev}(s_1)) \\
\theta\otimes\mathrm{id}\downarrow & & \downarrow\mathrm{id}\otimes\theta \\
\mathrm{ev}(s_3\star s_2)\otimes\mathrm{ev}(s_1) & & \mathrm{ev}(s_3)\otimes\mathrm{ev}(s_2\star s_1) \\
\theta\searrow & & \swarrow\theta \\
& \mathrm{ev}(s_3\star s_2\star s_1) &
\end{array}$$

*Left unitality.* Write $\mathbb{1}_* \coloneqq (*\mapsto\mathbb{1})$. We have equality at the bottom of the right diagram in $[*\times\mathcal{S}(a,b),\mathcal{V}]$, by left unitality for the Cat-enriched $\mathcal{S}_F$.

$$\begin{array}{ccc}
\mathbb{1}\boxtimes\mathbb{S}(a,b) & \xrightarrow{\left(j_b^{\mathcal{S}},\mathrm{id}\right)\boxtimes\mathrm{id}} & \mathbb{S}(b,b)\boxtimes\mathbb{S}(a,b) \\
\left(l^{\mathrm{Cat}},l^{\mathcal{V}}\right)\searrow & & \swarrow(\star,\theta) \\
& \mathbb{S}(a,b) &
\end{array}
\qquad
\begin{array}{ccc}
\mathbb{1}_*\boxtimes\mathrm{ev}^{a,b} & = & \left(\mathrm{ev}^{b,b}\circ j_b^{\mathcal{S}}\right)\boxtimes\mathrm{ev}^{a,b} \\
l^{\mathcal{V}}\Downarrow & & \Downarrow\theta\left(j^{\mathcal{S}}\times\mathrm{id}\right) \\
\mathrm{ev}^{a,b}\circ l^{\mathrm{Cat}} & = & \mathrm{ev}^{a,b}\left(j_b^{\mathcal{S}}(-)\star_{a,b,b}-\right)
\end{array}$$

Insert some $(\star,s)$. The triangle commutes by definition of $\theta$.

$$\begin{array}{ccc}
\mathbb{1}\otimes\mathrm{ev}(s) & = & \mathrm{ev}(\varnothing_b)\otimes\mathrm{ev}(s) \\
l^{\mathcal{V}}\downarrow & & \downarrow\theta \\
\mathrm{ev}(s) & = & \mathrm{ev}(\varnothing_b\star s)
\end{array}$$

*Right unitality.* Analogously to left unitality, unravel the diagram in $[\mathcal{S}(a,b)\times *,\mathcal{V}]$, use right unitality for $\mathcal{S}_F$ and insert some $(s,*)$. $\square$

# B. FUNCTORIALITY OF THE STRING CONSTRUCTION

Let colim be as in Proposition 3.3, let $\pi_*$ be as in Remark 4.12, and let $\mathrm{colim}_*$ be the change of base along the strong monoidal functor from Proposition 3.8. We have shown in Section 4 that the following diagram commutes pointwise for functors $F\colon \mathcal{I} \to \mathrm{Cat}_{\mathcal{V}}$. This section establishes functoriality of all constructions involved.

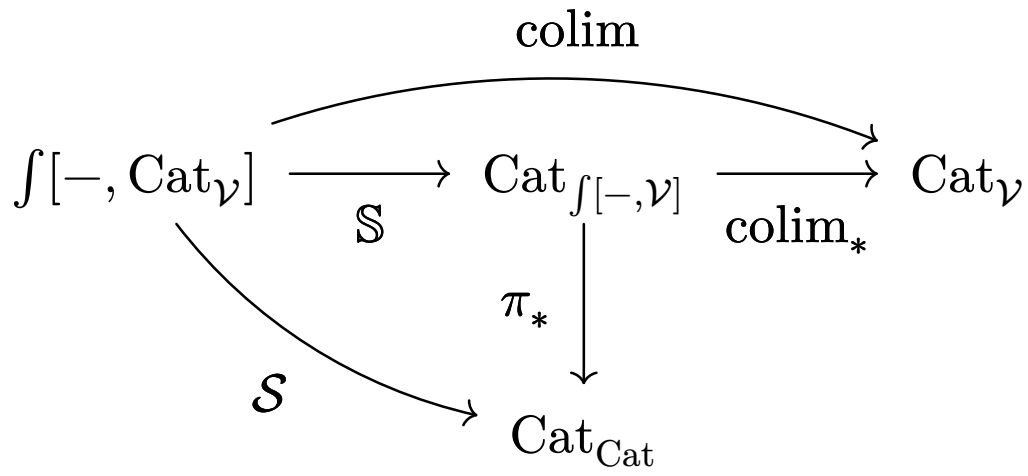

LEMMA B.1. Define the functor $\mathcal{S}\colon \int[-,\mathrm{Cat}_{\mathcal{V}}] \to \mathrm{Cat}_{\mathrm{Cat}}$ on objects by $\mathcal{S}(\mathcal{I}, F) \coloneqq \mathcal{S}_F$ as in Remark 4.6. Define the Cat-enriched functor $\mathcal{S}(H,\varphi)\colon \mathcal{S}_F \to \mathcal{S}_G$ on objects by $[x] \mapsto [\varphi_i(x)]$ and let the functor $\mathcal{S}(H,\varphi)_{a,b}$ between hom-categories be defined edgewise by $(i,x,y) \mapsto (H(i), \varphi_i(x), \varphi_i(y))$ and on empty strings by $\varnothing_a \mapsto \varnothing_{\mathcal{S}(H,\varphi)(a)}$.

PROOF. The assignment $\mathrm{colim}(\mathrm{Ob} \circ F) \to \mathrm{colim}(\mathrm{Ob} \circ G)$ is well-defined, since the composites

$$\mathrm{Ob}\, F(i) \xrightarrow{\mathrm{Ob}\,\varphi_i} \mathrm{Ob}\, G(Hi) \to \mathrm{colim}(\mathrm{Ob} \circ G)$$

define a cocone $\mathrm{Ob} \circ F \Rightarrow \Delta_{\mathcal{I}}\, \mathrm{colim}(\mathrm{Ob} \circ G)$. For a string $s \in \mathcal{S}_F(a,b)$, this also ensures coinciding endpoints in $\mathcal{S}(H,\varphi)_{a,b}(s)$. Let $\mathcal{S}(H,\varphi)_{a,b}$ send contractions to contractions and unit insertions to unit insertions in the evident way. Send a transport along $f\colon i \to j$ to the transport along $H(f)$, using naturality of $\varphi$.

$$\begin{array}{ccc} F(i) & \xrightarrow{\varphi_i} & G(H(i)) \\ {\scriptstyle F(f)}\downarrow & & \downarrow{\scriptstyle G(H(f))} \\ F(j) & \xrightarrow{\varphi_j} & G(H(j)) \end{array}$$

Composition and unit coherence of $\mathcal{S}(H,\varphi)$ follow straight from the definition. It is readily checked that $\mathcal{S}$ preserves compositions, since both $\mathcal{S}(K,\psi) \circ \mathcal{S}(H,\varphi)$ and $\mathcal{S}(K \circ H, \psi H \circ \varphi)$ send $(i,x,y)$ to the edge $(KHi, \psi_{Hi}(\varphi_i x), \psi_{Hi}(\varphi_i y))$. □

We now lift $\mathcal{S}$ along the change of enriching base $\pi_*\colon \mathrm{Cat}_{\int[-,\mathcal{V}]} \to \mathrm{Cat}_{\mathrm{Cat}}$ from Remark 4.12.

LEMMA B.2. Define the functor $\mathbb{S}\colon \int[-,\mathrm{Cat}_{\mathcal{V}}] \to \mathrm{Cat}_{\int[-,\mathcal{V}]}$ on objects by $\mathbb{S}(\mathcal{I},F) \coloneqq \mathbb{S}_F$ as in Definition 4.11. Abbreviating $\Phi \coloneqq \mathcal{S}(H,\varphi)$, define the $\int[-,\mathcal{V}]$-enriched functor $\mathbb{S}(H,\varphi)\colon \mathbb{S}_F \to \mathbb{S}_G$ by $\mathbb{S}(H,\varphi)_{a,b} \coloneqq \left(\Phi_{a,b}, \mathrm{ev}^{\varphi}_{a,b}\right)$ between hom-diagrams. The natural transformation

$$\mathrm{ev}^{\varphi}_{a,b}\colon \mathrm{ev}^{F}_{a,b} \Rightarrow \mathrm{ev}^{G}_{\Phi a,\Phi b} \circ \Phi_{a,b}$$

has components $\cdots \otimes (\varphi_i)_{x,y} \otimes \cdots$ for non-empty strings $s \coloneqq [\ldots,(i,x,y),\ldots]$, where $(\varphi_i)_{x,y}\colon F(i)(x,y) \to G(Hi)(\varphi_i x, \varphi_i y)$. If $a = b$, the $\varnothing_a$-component is $\mathrm{id}_{\mathbb{1}}$.

PROOF. $\mathrm{ev}^{\varphi}_{a,b}$ is natural in transports by naturality of $\varphi$, natural in unit insertions by unit coherence of the enriched functors $\varphi_i$, and natural in contractions by composition coherence of the enriched functors $\varphi_i$.

$$\begin{array}{ccc} F(i)(y,z) \otimes F(i)(x,y) & \xrightarrow{\circ} & F(i)(x,z) \\ {\scriptstyle \varphi_i \otimes \varphi_i}\downarrow & & \downarrow{\scriptstyle \varphi_i} \\ G(Hi)(\varphi_i y, \varphi_i z) \otimes G(Hi)(\varphi_i x, \varphi_i y) & \xrightarrow{\circ} & G(Hi)(\varphi_i x, \varphi_i z) \end{array}
\qquad
\begin{array}{ccc} & \mathbb{1} & \\ {\scriptstyle j^{Fi}_x}\swarrow & & \searrow{\scriptstyle j^{GHi}_{\varphi_i x}} \\ F(i)(x,x) & \xrightarrow{\varphi_i} & G(Hi)(\varphi_i x, \varphi_i x) \end{array}$$

We verify composition coherence of $\mathbb{S}(H,\varphi)$ by unraveling the relevant compositions in $\int[-,\mathcal{V}]$ and using composition coherence of $\Phi$ to obtain a diagram in $[\mathcal{S}_F(b,c)\times\mathcal{S}_F(a,b),\mathcal{V}]$. Inserting strings $s'$ and $s$, we obtain the following commutative diagram in $\mathcal{V}$, where we suppress superscript endpoints, and $\theta^F$ and $\theta^G$ come from Lemma 4.10.

$$\begin{array}{ccccc}
\mathrm{ev}^F(s')\otimes\mathrm{ev}^F(s) & & \xrightarrow{\theta^F} & & \mathrm{ev}^F(s'\star s)\\
\downarrow{\scriptstyle \mathrm{ev}^\varphi\otimes\mathrm{ev}^\varphi} & & & & \downarrow{\scriptstyle \mathrm{ev}^\varphi}\\
\mathrm{ev}^G(\Phi s')\otimes\mathrm{ev}^G(\Phi s) & \xrightarrow{\theta^G} & \mathrm{ev}^G(\Phi s'\star\Phi s) & = & \mathrm{ev}^G\Phi(s'\star s)
\end{array}$$

For unit coherence of $\mathbb{S}(H,\varphi)$, unravel the diagram in $\int[-,\mathcal{V}]$ to obtain a diagram in $[*,\mathcal{V}]$, which commutes, since $\mathrm{ev}^G(\varnothing_{\Phi a})=\mathrm{ev}^F(\varnothing_a)=\mathbb{1}$. ☐

Proposition B.3. By Proposition 3.3, a choice of colimit per diagram defines a functor $\mathrm{colim}\colon\int[-,\mathrm{Cat}_\mathcal{V}]\to\mathrm{Cat}_\mathcal{V}$. Choosing $\mathrm{colim}(F):=\mathrm{colim}_*(\mathbb{S}_F)$, this functor agrees with $\mathrm{colim}_*\circ\mathbb{S}$ on morphisms.

Proof. Let $\iota^F$ and $\iota^G$ be the cocones of $\mathrm{colim}_\mathcal{I}(F)$ and $\mathrm{colim}_\mathcal{J}(G)$ as in Lemma 4.15 respectively, and write $\Phi:=\mathcal{S}(H,\varphi)$ and $\tilde\Phi:=\mathbb{S}(H,\varphi)$. Then Proposition 3.3 defines $\mathrm{colim}(H,\varphi)\colon\mathrm{colim}(F)\to\mathrm{colim}(G)$ to be the unique morphism from Lemma 4.16, satisfying commutativity of the left diagram in $[\mathcal{I},\mathrm{Cat}_\mathcal{V}]$, spelled out component-wise in the right diagram in $\mathrm{Cat}_\mathcal{V}$.

$$\begin{array}{ccc}
F & \overset{\varphi}{\Longrightarrow} & G\circ H\\
\Downarrow{\scriptstyle\iota^F} & & \Downarrow{\scriptstyle\iota^G H}\\
\Delta\,\mathrm{colim}\,F & \overset{\Delta\,\mathrm{colim}_*\tilde\Phi}{\Longrightarrow} & \Delta\,\mathrm{colim}\,G
\end{array}
\qquad
\begin{array}{ccc}
F(i) & \xrightarrow{\varphi_i} & G(H(i))\\
\downarrow{\scriptstyle\iota^F_i} & & \downarrow{\scriptstyle\iota^G_{H(i)}}\\
\mathrm{colim}_*\mathbb{S}_F & \xrightarrow{\mathrm{colim}_*\tilde\Phi} & \mathrm{colim}_*\mathbb{S}_G
\end{array}
\tag{4}$$

Now $\mathrm{colim}_*\tilde\Phi\colon\mathrm{colim}_*\mathbb{S}_F\to\mathrm{colim}_*\mathbb{S}_G$ sends $\mathrm{colim}(\mathrm{Ob}\circ F)\to\mathrm{colim}(\mathrm{Ob}\circ G)$ as in Lemma B.1. Let $\kappa^F$ and $\kappa^G$ denote the cocones of $\mathrm{ev}^F_{a,b}$ and $\mathrm{ev}^G_{\Phi a,\Phi b}$ respectively. Recall from Lemma 2.5 that $\big(\mathrm{colim}_*\tilde\Phi\big)_{a,b}=\mathrm{colim}\,\tilde\Phi_{a,b}=\mathrm{colim}\big(\Phi_{a,b},\mathrm{ev}^\varphi_{a,b}\big)$, i.e. the unique morphism satisfying commutativity of the left diagram in $[\mathcal{S}_F(a,b),\mathcal{V}]$. Inserting some string $s:=[(i,x,y)]$, we obtain the right diagram.

$$\begin{array}{ccc}
\mathrm{ev}^F_{a,b} & \overset{\mathrm{ev}^\varphi_{a,b}}{\Longrightarrow} & \mathrm{ev}^G_{\Phi a,\Phi b}\circ\Phi_{a,b}\\
\Downarrow{\scriptstyle\kappa^F} & & \Downarrow{\scriptstyle\kappa^G\Phi_{a,b}}\\
\Delta\,\mathrm{colim}\,\mathbb{S}_F(a,b) & \overset{\Delta\,\mathrm{colim}\,\tilde\Phi_{a,b}}{\Longrightarrow} & \Delta\,\mathrm{colim}\,\mathbb{S}_G(\Phi a,\Phi b)
\end{array}
\qquad
\begin{array}{ccc}
F(i)(x,y) & \xrightarrow{(\varphi_i)_{x,y}} & G(Hi)(\varphi_i x,\varphi_i y)\\
\downarrow{\scriptstyle\kappa^F_s} & & \downarrow{\scriptstyle\kappa^G_{\Phi(s)}}\\
\mathrm{colim}\,\mathbb{S}_F(a,b) & \xrightarrow{\mathrm{colim}\,\tilde\Phi_{a,b}} & \mathrm{colim}\,\mathbb{S}_G(\Phi a,\Phi b)
\end{array}$$

Evidently, the right diagram concludes commutativity of Equation 4. ☐

# C. Alternative rigidification of simplices

Let Grph denote the category of reflexive graphs. Let $F\colon\mathrm{Grph}\to\mathrm{Cat}$ be the free and $U\colon\mathrm{Cat}\to\mathrm{Grph}$ the forgetful functor. The free-forgetful adjunction $F\dashv U$ gives rise to the comonad resolution $(FU)_\bullet(\mathcal{C})$ given by $[n]\mapsto(FU)^{n+1}(\mathcal{C})$, defining a simplicially enriched category.

$$FU\mathcal{C}\;\overset{\varepsilon FU,\;F\eta U,\;FU\varepsilon}{\rightleftarrows}\;FUFU\mathcal{C}\;\overset{\varepsilon FUFU,\;F\eta UFU,\;FU\varepsilon FU,\;FUF\eta U,\;FUFU\varepsilon}{\rightleftarrows}\;FUFUFU\mathcal{C}\;\cdots$$

Remark C.1. While the comonad resolution $(FU)_\bullet[n]$ provides insightful intuition, Lurie's model for $\mathfrak{C}\Delta^n$ in Definition 6.7 is more suitable for computations.

Example C.2. We compute the hom-sSet $((FU)_\bullet[3])(0,3)$, beginning with $[3]$ on the left. Applying the forgetful functor $U\colon \mathsf{Cat} \to \mathsf{Grph}$, we forget that $hg$, $hgf$ and $gf$ are compositions and name the resulting arrows $(hg)$, $(hgf)$ and $(gf)$ respectively. Applying the free functor $F\colon \mathsf{Grph} \to \mathsf{Cat}$, we add new compositions, depicted on the right.

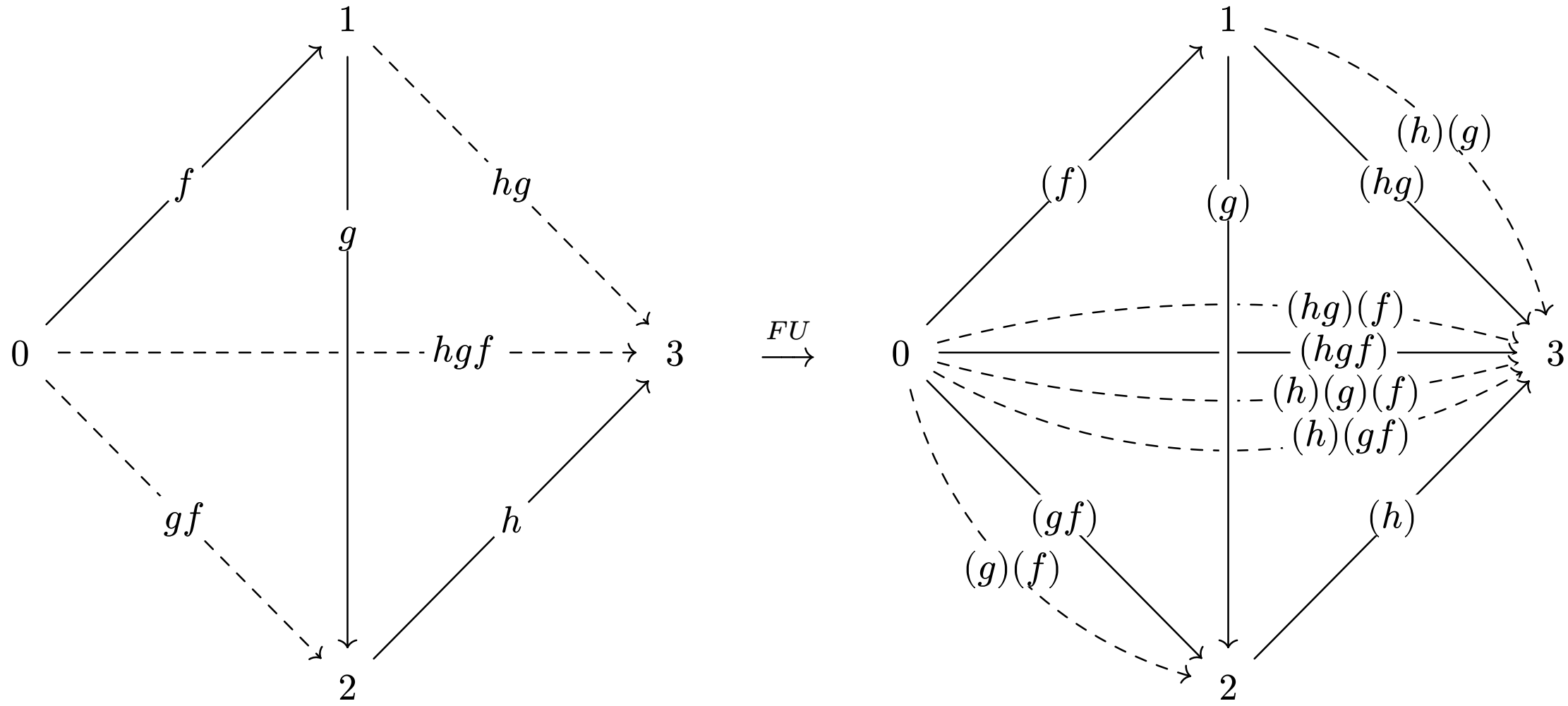

The face maps $(FU)^i\varepsilon(FU)^j$ remove the parentheses contained in exactly $i$ others. The degeneracy maps $F(UF)^k\eta(UF)^jU$ double the parentheses contained in exactly $k$ others. In general, $((FU)_\bullet[n])(i,j)$ is empty if $i > j$, a point if $i = j$ and a cube $(\Delta^1)^{j-i-1}$ if $i < j$. Thus, $((FU)_\bullet[3])(0,3) = \Delta^1 \times \Delta^1$.

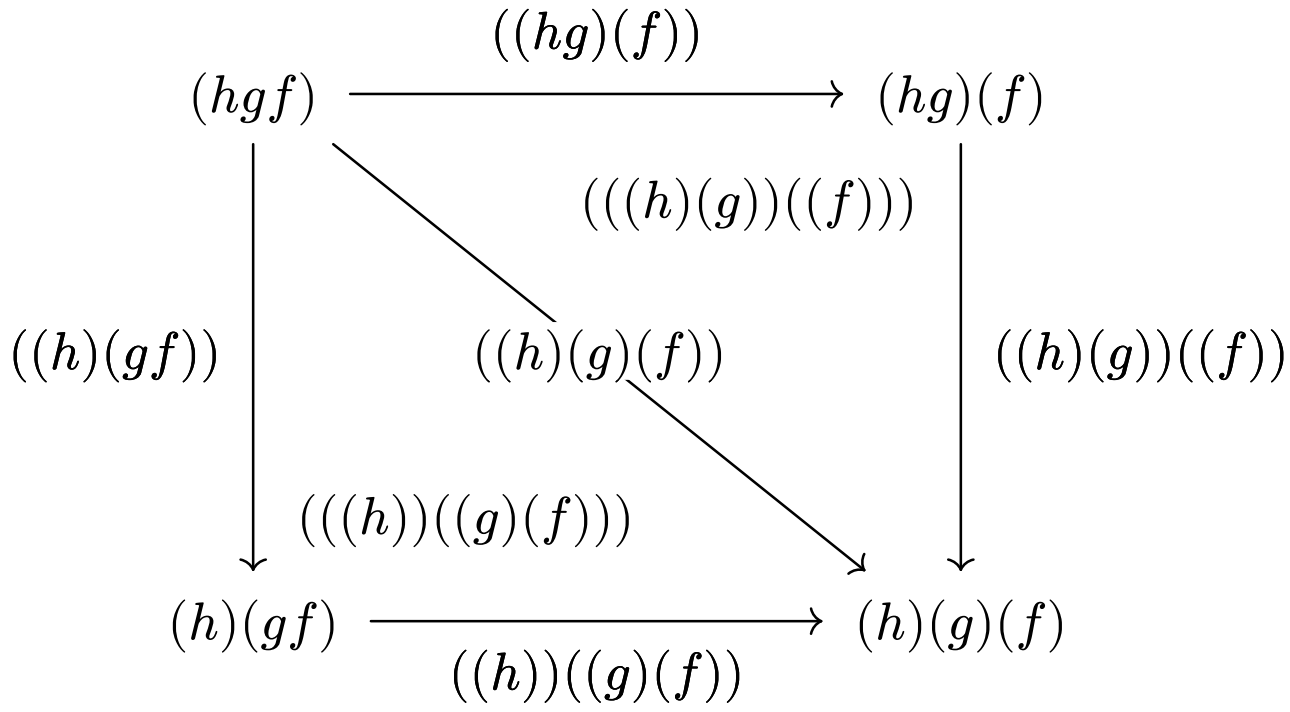

> Lemma C.3 [DS11, Lemma 2.5]. For each $n$, the comonad resolution $(FU)_\bullet[n]$ is isomorphic to the rigidification $\mathfrak{C}\Delta^n$ from Definition 6.7.

Example C.4. From Example C.2, identify the 0-simplex $0 \xrightarrow{(gf)} 2 \xrightarrow{(h)} 3$ with $\{0,2,3\} \in (NP(0,3))_0$. Identify the 2-simplex $(((h))((g)(f)))$ with $\{0,3\} \subseteq \{0,2,3\} \subseteq \{0,1,2,3\}$.

> Lemma C.5. Let $\Delta$ denote the Yoneda embedding $\Delta \to [\Delta^{\mathrm{op}}, \mathsf{Set}]$ given by $[n] \mapsto \Delta^n$. The isomorphism from the proof of Lemma C.3 (see [DS11, §A.5]) extends to a natural isomorphism $(FU)_\bullet \simeq \mathfrak{C} \circ \Delta$ of functors $\Delta \to \mathsf{Cat}_{\mathsf{sSet}}$.

Proof. We reverse the order of flags to avoid renaming the sets $I_r$. Let $f\colon [m] \to [n]$ in $\Delta$. Let the components of $\alpha\colon (FU)_\bullet \Rightarrow \mathfrak{C} \circ \Delta$ be defined as in Lemma C.3.

$$\begin{array}{ccc} ((FU)_\bullet[m])(i,j) & \xrightarrow{(FU)_\bullet(f)_{i,j}} & ((FU)_\bullet[n])(fi,fj) \\ \Big\downarrow{\scriptstyle \alpha_{[m]}} & & \Big\downarrow{\scriptstyle \alpha_{[n]}} \\ NP(i,j) & \xrightarrow{\left(N_*\tilde{f}\right)_{i,j}} & NP(fi,fj) \end{array}$$

*Induction hypothesis.* For a $k$-simplex $g \in \big((FU)^{k+1}[m]\big)(i,j)$, the following diagram commutes.

$$\begin{array}{ccc} g & \xmapsto{(FU)^{k+1} f} & \big((FU)^{k+1} f\big)(g) \\ \alpha_{[m]} \downarrow & & \downarrow \alpha_{[n]} \\ (I_k \subseteq \cdots \subseteq I_0) & \xmapsto{N_* \tilde{f}} & (f(I_k) \subseteq \cdots \subseteq f(I_0)) \end{array}$$

*Base case.* A morphism $g \in [m](i, j)$ can be regarded as a "$(-1)$-simplex" sent to the "empty flag".

*Induction step.* A $k$-simplex $g \in \big((FU)^{k+1}[m]\big)(i, j)$ is a free composition of edges $(g_r) \in \big(U(FU)^k\big)[m]$

$$i \xrightarrow{(g_0)} i_1 \xrightarrow{(g_1)} i_2 \to \cdots \to i_\ell \xrightarrow{(g_\ell)} j$$

where $g_0, ..., g_\ell$ are morphisms in $(FU)^k[m]$. Say $\alpha_{[m]}$ sends the $(k-1)$-simplex $g_\ell \circ \cdots \circ g_0$ to the flag $(I_{k-1} \subseteq \cdots \subseteq I_0)$. Let $I_k := \{i, i_1, ..., i_\ell, j\} \in P(i, j)$ and define $\alpha_{[m]}(g) := (I_k \subseteq I_{k-1} \subseteq \cdots \subseteq I_0) \in NP(i, j)_k$. Now $N_* \tilde{f}$ sends this flag to $(f(I_k) \subseteq \cdots \subseteq f(I_0))$. Let $h_r := \big((FU)^k f\big)(g_r)$ and note that $\big((FU)^{k+1} f\big)((g_r)) = (h_r)$. Now $(FU)^{k+1}(f)$ sends $g$ to the free composition

$$f(i) \xrightarrow{(h_0)} f(i_1) \to \cdots \to f(i_\ell) \xrightarrow{(h_\ell)} f(j).$$

Applying the induction hypothesis to the $(k-1)$-simplex $g_\ell \circ \cdots \circ g_0$, we get

$$\alpha_{[n]}(h_\ell \circ \cdots \circ h_0) = \alpha_{[n]}\Big(\big((FU)^k f\big)(g_\ell \circ \cdots \circ g_0)\Big) = (f(I_{k-1}) \subseteq \cdots \subseteq f(I_0)).$$

Since $f(I_k) = \{f(i), f(i_1), ..., f(i_\ell), f(j)\}$, this concludes the proof. □